\documentclass{amsart}

\usepackage{latexsym}
\usepackage{amsmath}
\usepackage{dsfont}
\usepackage[all]{xy}

\newtheorem{teo}{Theorem}[section]
\newtheorem{lem}[teo]{Lemma}
\newtheorem{corol}[teo]{Corollary}
\newtheorem{defi}[teo]{Definition}
\newtheorem{exa}[teo]{Example}
\newtheorem{prop}[teo]{Proposition}
\newtheorem{rem}[teo]{Remark}

\newtheorem{notation}[teo]{Notation}

\newcommand{\beq}{\begin{equation}}
\newcommand{\eeq}{\end{equation}}

\newcommand\e{{\rm e}}
\newcommand{\A}{{\mathsf{A}}}

\newcommand{\T}{{\mathcal{T}}}
\renewcommand{\P}{{\mathcal{P}}}
\newcommand{\ec}{{\mathsf e}}
\renewcommand{\ge}{{\mathsf g}}

\newcommand{\ord}{{\mathsf o}}
\renewcommand{\det}{{\rm det}}

\newcommand{\N}{{\mathds{N}}}
\newcommand{\Z}{{\mathds{Z}}}

\newcommand{\C}{{\mathds{C}}}

\newcommand{\K}{{\mathds{K}}}

\newcommand{\Rs}[2]{{{\rm Res}_{#1}}_{_{_{_{\hspace{-.23in}\scriptscriptstyle{s=#2}}}}}\hspace{0.06in}}

\DeclareMathOperator*{\Rz}{Res_0}

\DeclareMathOperator*{\Ru}{Res_1}

\DeclareMathOperator*{\Rd}{Res_2}

\newcommand{\Res}[2]{{{\rm Res}_{#1}}_{_{\hspace{-.35in}\begin{array}{c}\scriptstyle{s=#2}\end{array}}}\hspace{-.007in}}
\renewcommand{\Re}{{\rm Re}}

\title{Zeta determinant for double sequences of spectral type}

\author{M. Spreafico}
\address{ICMC, Universidade S\~{a}o Paulo, S\~{a}o Carlos, Brazil, {\tt mauros@icmc.usp.br}.}

\begin{document}

\maketitle

\begin{abstract}
We study the spectral functions, and in particular the zeta function, associated to a class of sequences of complex numbers, called of spectral type. We investigate the decomposability of the zeta function associated to a double sequence with respect to some simple sequence, and we provide a technique for obtaining the first terms in the Laurent expansion at zero of the zeta function associated to a double sequence. We particularize this technique to the case of sums of sequences of spectral type, and we give two applications: the first concerning some special functions appearing in number theory, and the second the functional determinant of the Laplace operator on a product space.  
\end{abstract}






\section{Introduction}
\label{s0}

Let $S=\{\lambda\}$ be a sequence of complex numbers. Under some mild requirements on $S$, it is possible to define the zeta function associated to $S$ by 
\beq\label{z}
\zeta(s,S)=\sum_{\lambda\in S}^\infty \lambda^{-s},
\eeq
where $s$ is a complex variable, and the series in equation (\ref{z}) converges in some domain of the complex plane. For example, if $S=\N_0=\{1,2,\dots\}$ are the positive integers, then $\zeta(s,S)=\zeta_R(s)$ is the Riemann zeta function, if $S$ are the norms of the integral ideals of a number field $\K$, then $\zeta(s,S)=\zeta_\K(s)$ is the Dedekind zeta function of $\K$; in both cases the series in equation (\ref{z}) converges if $\Re(s)>1$. If $S$ is the spectrum of the Laplace operator on a compact Riemannian manifold of dimension $m$, then the series in equation (\ref{z}) converges if $\Re(s)>\frac{m}{2}$.

We introduced in \cite{Spr4} a class of abstract sequences of complex numbers, that we called of spectral type, and we  investigated the properties of the spectral functions associated to the sequences of that class.  In particular, we studied the expansion of the zeta function at zero. In fact, it is well know that the determination of the value of the derivative of the zeta function at zero is a major problem  in various areas of mathematics and theoretical physics. In number theory, it is related with the ideal class number, while in analysis and geometry it is used in order to define the functional determinant of a linear operator. In a series of related works \cite{Spr2, Spr3, Spr5, Spr7, SZ}, we studied various different applications of the approach introduced in \cite{Spr4}. In particular, in \cite{Spr3} and \cite{Spr5} some double sequences appeared, and we introduced some techniques that allowed to tackle the expansion at zero of the associated zeta functions. These techniques permitted to  obtain very useful information in the different problems analyzed in those works. This motivates the investigation of a general purpose technique, that permits to deal with a large class of double or even multiple sequences. This is the aim of the present work, that is structured as follows. In Section \ref{s1}, we review the definition of sequences of spectral type given in \cite{Spr4}, and we give some further properties of the zeta function associated to the sequences of this type. In Section \ref{s2}, we consider double (or multiple) sequences. We define a criterium for decomposing a double sequence with respect to some simple sequence of spectral type (Definition \ref{spdec}), and we prove our main result in Theorem \ref{sdl}, where we give the expansion at zero of the zeta function associated to a decomposable double  sequence. 
In Section \ref{s3}, we apply Theorem \ref{sdl} to investigate the case of sums of sequences of spectral type, namely zeta functions of the form $\zeta(s,S_{(1)}+S_{(2)})=\sum_{m,n=1}^\infty (\lambda_{1,m}+\lambda_{2,n})^{-s}$, where $S_{(j)}=\{\lambda_{j,m}\}$ are two sequences of spectral type. We obtain in Theorem \ref{tt1} the expansion of the zeta function $\zeta(s,S_{(1)}+S_{(2)})$ at zero. Beside the natural interest for  studying this problem in zeta regularization theory \cite{CQ, JL1, HQS}, motivations arise from different areas of mathematics and theoretical physics.  We discuss two applications in the last Section \ref{s4}, the first concerning some special functions appearing in number theory \cite{Zag1, Spr5}, the second the zeta determinant of the Laplace operator on a product manifold \cite{FdG, OS}.

\section{Zeta invariants for sequences of spectral type}
\label{s1}

In this section, we recall from \cite{Spr4} (see also \cite{Vor, CQ, HQS, JL1} for other approach to zeta regularized products) the definition of sequences of spectral type, and we generalize some results of that work. Let $S=\{\lambda_n\}_{n=1}^\infty$ be a sequence
of non vanishing complex numbers with unique accumulation point at infinity. We ordered $S$ as $0<|\lambda_1|\leq |\lambda_2|\leq\dots$ (if we need to include the number zero, we use the notation $\lambda_0=0$). 
The positive number
\[
\ec=\ec(S)={\rm limsup}_{n\to\infty} \frac{\log n}{\log |\lambda_n|},
\]
is called the {\it exponent of convergence} of $S$. We are only interested in sequences with $\ec(S)<\infty$. If this is the case, then there exists a least integer $p$ such that the series $\sum_{n=1}^\infty \lambda_n^{-p-1}$ converges absolutely.  We assume $\ec-1<p\leq \ec$, and hence we have $p=[\ec]$ (the integer part). The integer $p$ is called the {\it genus} of the sequence $S$, and we use the notation $\ge=\ge(S)=p$. The Weierstrass canonical product 
\[
p(z,S)=\prod_{n=1}^\infty\left(1+\frac{z}{\lambda_n}\right)\e^{\sum_{j=1}^{\ge(S)}\frac{(-1)^j}{j}\frac{z^j}{\lambda_n^j}},
\]
converges uniformely and absolutely in any bounded closed region of the complex plane, and is an integral function of order $\ge(S)$ which vanishes if and only if $z=-\lambda_n$ for some $n$. We call the open subset $\rho(S)=\C-S$ of the complex plane the  {\it resolvent set} of $S$. We define the {\it Gamma function} associated to $S$  by the canonical product $\Gamma(z,S)=\frac{1}{p(z,S)}$.
For further convenience, we  use the variable $\lambda=-z$ for the Gamma function, with the convention that the variable $-\lambda$ is purely real on the negative part of the real axes. When necessary, in order to define the branch of analicity of a meromorphic function, the domain for $\lambda$ will be some  open subset of  $\C-[0,\infty)$ in the complex plane. With this notation,  for all $\lambda\in\rho(S)$, the Gamma functions is
\beq\label{gamma}
\frac{1}{\Gamma(-\lambda,S)}=\prod_{n=1}^\infty\left(1+\frac{-\lambda}{\lambda_n}\right)\e^{\sum_{j=1}^{\ge(S)}\frac{(-1)^j}{j}\frac{(-\lambda)^j}{\lambda_n^j}}.
\eeq

We introduce now the definition of sequences of spectral type. Conditions for a sequence to be of spectral type can be found in \cite{JL1} or \cite{Spr4}.

\begin{defi}\label{sss} Let $S=\{\lambda_n\}_{n=1}^\infty$ be a sequence of non vanishing complex numbers with unique accumulation point at infinity and exponent of convergence $\ec(S)<\infty$. We call $S$ a {\it sequence of spectral type ($S$-type)}  if the following conditions hold:
\begin{enumerate} 
\item there exist $c>0$ and $0<\theta<\pi$, such that the sequence $S$ is contained in the interior of the sector $\Sigma_{\theta,c}=\left\{z\in \C~|~|\arg(z-c)|\leq \frac{\theta}{2}\right\}$;
\item the logarithm of the associated Gamma function has an 
asymptotic expansion for large $\lambda\in AD_{\theta,c}(S)=\C-\Sigma_{\theta,c}$ of the following form
\[
\log\Gamma(-\lambda,S)=\sum_{\alpha} \sum_{k=0}^{K_\alpha} a_{\alpha,k}(-\lambda)^\alpha \log^k (-\lambda) +o((-\lambda)^{\alpha_N}),
\]
where $\{\alpha\}$ is a decreasing sequence of real numbers $\alpha_0> \alpha_1>\dots >\alpha_N\geq-\infty$, and $k=0,1,\dots,K_\alpha\in \N$, for each $\alpha$.
\end{enumerate}

If $N$ is finite, we call the number $\alpha_N$ the {\it order} of the sequence $S$, and we use the notation $\ord(S)$. We say that $S$ has {\it infinite order} if $N$ is not finite, and we write $\ord(S)=-\infty$. We call the open set $AD_{\theta,c}(S)$ the {\it asymptotic domain} of $S$. 
\end{defi}

\begin{rem}\label{r1} The point $\lambda=0$ belongs by definition to the domain of analicity of $\Gamma^{-1}(-\lambda,S)$ for a sequence of spectral type, and $\log\Gamma(-\lambda,S)$ has a zero of order $\ge(S)+1$ at $\lambda=0$.
\end{rem}

We define some further spectral functions associated to a sequence $S$ of spectral type. The {\it zeta function}, defined by the uniformely convergent series 
\[
\zeta(s,S)=\sum_{n=1}^\infty \lambda_n^{-s},
\]
for $\Re(s)>\ec(S)$, and by analytic continuation elsewhere, the {\it heat function} ($t>0$) 
\beq\label{heat}
f(t,S)=1+\sum_{n=1}^\infty \e^{-t\lambda_n},
\eeq
and, for each non negative integer $k$, and $\lambda\in\rho(S)$, the following functions
\begin{align*}
r_k(\lambda,S)=&k!\sum_{n=1}^\infty\left(\frac{(-1)^k}{(\lambda-\lambda_n)^{k+1}}+
\sum_{j=0}^{\ge(S)-k-1} \binom{j+k}{j}\frac{\lambda^j}{\lambda^{j+k+1}_n}\right).
\end{align*}

These spectral functions are clearly strictly related to each others (under
some general regularity's assumptions). In particular the following complex integral representations hold (see \cite{Spr2} and \cite{Spr4} Proposition 2.4). Here we use the notation $\Sigma_{\theta,c}=\left\{z\in \C~|~|\arg(z-c)|\leq \frac{\theta}{2}\right\}$, 
with $c\geq \delta> 0$, $0< \theta<\pi$. We use $AD_{\theta,c}=\C-\Sigma_{\theta,c}$, for the complementary (open) domain and $\Lambda_{\theta,c}=\partial \Sigma_{\theta,c}=\left\{z\in \C~|~|\arg(z-c)|= \frac{\theta}{2}\right\}$, oriented counter clockwise, for the boundary.

\begin{lem} \label{f} For $\Re(s)>\ec(S)$,
\begin{align}
\label{f0}f(t,S)-1&=\frac{1}{2\pi i}\int_{\Lambda_{\theta,c}}\e^{-\lambda t}
r_0(\lambda,S) d\lambda\\
\label{f0b}&=-\frac{t}{2\pi
i}\int_{\Lambda_{\theta,c}}\e^{-\lambda t} \log \Gamma(-\lambda,S)
d\lambda,\\
\label{f1}\zeta(s,S)&=\frac{1}{\Gamma(s)}\int_0^\infty t^{s-1}(f(t,S)-1)dt\\
\label{f1b}&=\frac{s}{\Gamma(s)}\int_0^\infty t^{s-1}\frac{1}{2\pi
i}\int_{\Lambda_{\theta,c}}\frac{\e^{-\lambda t}}{-\lambda} \log
\Gamma(-\lambda,S) d\lambda dt.
\end{align}
\end{lem}

\begin{lem} \label{l0} For all $k$, 
\[
r_k(\lambda,S)=-\frac{d^{k+1}}{d\lambda^{k+1}}\log\Gamma(-\lambda,S).
\]
\end{lem}

\begin{lem}\label{r2} If $k=\ge=\ge(S)$, then
\[
r_{\ge}(\lambda,S)=-\ge!\zeta(\ge+1,S-\lambda) =- \ge!\sum_{n=1}^\infty
(\lambda_n-\lambda)^{-\ge-1}.
\]

This is a uniformly convergent series and we can take the limit
for $\lambda\to \infty$ that is 0. This implies that
$r(\lambda,S)$ can not have a term like $(-\lambda)^m$ and
consequently $\log\Gamma(-\lambda,S)$ can not have a term like
$(-\lambda)^{m+1}$, with $m\geq \ge$.
\end{lem}

\begin{lem}\label{vecchio} Let $S$ be a
sequence with finite exponent contained in a positive sector $\Sigma_{\theta,c>0}$. Then
the heat function associated to $S$ has an asymptotic expansion of
order $\delta'\geq-\infty$ with respect to the asymptotic sequence (see for example \cite{Mur} for definition) 
$\{t^{-\delta}\log^l t\}_\delta$ as $t\to 0^+$ if and only if the
logarithmic $\Gamma$-function associated to $S$ has an asymptotic
expansion of order $\alpha'\geq -\infty$ with respect to the
asymptotic sequence $\{(-\lambda)^{\alpha}\log^k
(-\lambda)\}_\alpha$ as $\lambda\to \infty$ uniformly in $\lambda$
for $\lambda\in B$, where $B$ is any unbounded region contained in
the complement of $\Sigma_{\theta,c}$ (and $\delta'=-\infty$ if
and only if $\alpha'=-\infty$). The order $\delta'$
and the coefficients in the expansion of the heat function can be
obtained from $\alpha'$ and the coefficients in the
expansion of the logarithmic Gamma function, but not viceversa. 
\end{lem}

Even if many useful results hold for a generic sequence of $S$-type, we consider here only sequence of a particular subclass as defined below (see also \cite{BS2} Section 9).

\begin{defi}\label{reg} A sequence of spectral type $S$ is called {\it regular} if the coefficients $a_{\alpha,k}$ in the expansion of $\log\Gamma(-\lambda,S)$ 
vanish for all $k\not=0,1$.
\end{defi}

\begin{rem} \label{r5} Let $S$ be a regular sequence of spectral type with $\ord(S)=\alpha_N$. Then, 
\[
\log\Gamma(-\lambda,S)=\sum_{j=0}^N a_{\alpha_j,0}(-\lambda)^{\alpha_j}+\sum_{j=0}^N a_{\alpha_j,1}(-\lambda)^{\alpha_j}\log(-\lambda)+O((-\lambda)^{\alpha_N}),
\]
for large $\lambda$ in $AD_{\theta,a}(S)$, and where $\alpha_0>\alpha_1>\dots>\alpha_N=\ord(S)$.
\end{rem}

\begin{rem}\label{l2} If $S$ is a regular sequence of $S$-type, then $\alpha_0\leq \ec(S)$, and $\alpha_0<\ge(S)+1$.
\end{rem}

We give now the main results on regular sequences of $S$-type in this context. The first result follows from equation (\ref{f1}) and properties of special functions.

\begin{prop}  \cite{Spr4} Let $S$ be a regular sequence of spectral type with finite
exponent and order $\alpha_N$. Then, the associated heat function has the following asymptotic expansion
for $t\to 0^+$
\[
f(t,S)-1=\sum_{j=0}^N \sum_{k=0}^1
c_{\alpha_j,k}t^{-\alpha_j}\log^k t +o(t^{-\alpha_N}),
\]
where: $c_{\alpha_j,0}=\frac{1}{\Gamma(-\alpha_j)}\left(a_{\alpha_j,0}+\psi(-\alpha_j)a_{\alpha_j,1}\right)$, $c_{\alpha_j,1}=-\frac{a_{\alpha_j,1}}{\Gamma(-\alpha_j)}$, and $\psi$ is the Digam\-ma function.
%
\label{coeff}\end{prop}

Note that the equations among the coefficients given in 
Proposition \ref{coeff} are invertible provided $\alpha_j\not\in
\N$; moreover, $c_{\alpha_j,1}=0$, whenever $\alpha_j\in
\N$. Note also that, the formulas in Proposition
\ref{coeff} for the coefficients only hold for regular sequences
of spectral type. For a generic $S$-sequence, it is not true for
example that the coefficients $c_{\alpha_j,1}$ are null for
$\alpha_j$ a natural number.

Next, we give formulas for the  zeta invariants
of a regular sequences. This result corresponds to the classical result for the functional zeta function of a differential operator. See also Theorem 3.4, 3.6 of \cite{JL1} for a further generalization to a larger class of sequences. 

\begin{prop} Let $S$ be a regular sequence of spectral type of order
$\alpha_N$. Then, the associated
zeta function is regular in the complex half-plane
$\Re(s)>\alpha_N-\epsilon$ (positive small $\epsilon$), up to a
finite set of poles. The poles in the half plane
$\Re(s)>\alpha_N-\epsilon$ are at most $N+1$, are located at
$s=\alpha_N,\alpha_{N-1},\dots, \alpha_0\leq \ec(S)$, and are of
order at most $2$ with residues
\begin{align*}
\Rd_{s=\alpha_j}\zeta(s,S)&=-\frac{c_{\alpha_j,1}}{\Gamma(\alpha_j)}
\hspace{.3in} \left(=0, ~{\rm if}~\alpha_j\in\Z\right),\\
\Ru_{s=\alpha_j}\zeta(s,S)&=\left\{\begin{array}{ll}\frac{c_{\alpha_j,0}}{\Gamma(\alpha_j)}
+c_{\alpha_j,1}\frac{\psi(\alpha_j)}{\Gamma(\alpha_j)},&
~\alpha_j\not= 0,-1,-2,\dots,-[|\alpha_N|],\\
(-1)^{-\alpha_j+1}(-\alpha_j)! c_{\alpha_j,1},& ~\alpha_j=
-1,-2,\dots,-[|\alpha_N|], \end{array}\right.\\
\Rz_{s=\alpha_j}\zeta(s,S)&=(-1)^{\alpha_j}(-\alpha_j)!c_{\alpha_j,0},\hspace{.6in}
\alpha_j=-1,-2,\dots,-[|\alpha_N|];
\end{align*}
in particular, if $\alpha_N\leq 0$, then $s=0$ is a regular point with $\zeta(0,S)=c_{0,0}$
\label{zeta1}
\end{prop}

It is important to observe that all the formulas given Proposition \ref{zeta1} can be written using exclusively 
the coefficients $a_{\alpha_j,k}$ appearing in the asymptotic
expansion of the $\Gamma$-function. This follows from Proposition
\ref{coeff} or by direct computations (see also Proposition \ref{aaab3} below).

\begin{teo} \cite{Spr4} If $S$ is a regular sequence of spectral type of order
$\alpha_N\leq 0$, then the associated zeta function 
is regular at $s=0$, and near $s=0$
\[
\zeta(s,S)=-a_{0,1}-a_{0,0}s+O(s^2).
\]
\label{zeta2}
\end{teo}

We call a regular sequence of spectral type {\it simply regular}
if the unique logarithmic terms appearing in the expansion of
$\log\Gamma(-\lambda,S)$ are of the form
$(-\lambda)^k\log(-\lambda)$, with integer $k$, and {\it totally
regular} if also $k\geq 0$. It is clear that the zeta function of
a simply regular sequence of order $\alpha_N$ has at most simple
poles for $\Re(s)\geq \alpha_N$. It is also easy to see that, while there
can be logarithmic terms in the expansion of the heat function of
a simply regular sequence, namely of the type $t^k \log t$, with
negative integer $k$, we have $c_{\alpha_j,1}=0$ for all
$\alpha_j$ for a totally regular sequence.

\begin{rem} \label{form} When $S$ is a totally regular sequence of spectral type with order
$\alpha_N$ and genus $\ge$, the asymptotic expansion of
the associated spectral functions reads:
\begin{align*}
\log\Gamma(-\lambda,S)&=\sum_{j=0}^\ge a_{j,1}(-\lambda)^j \log
(-\lambda) +\sum_{j=0}^N
a_{\alpha_j,0}(-\lambda)^{\alpha_j}+o((-\lambda)^{\alpha_N}),\\
f(t,S)-1&=\sum_{j=0}^N c_{\alpha_j,0}t^{-\alpha_j}+o(t^{-\alpha_N}).
\end{align*}
\end{rem}

The main results concerning totally regular sequences are the
following ones. First,  we have formulas that express the
coefficients $a_{\alpha_j,k}$ as functions of the coefficients
$c_{\alpha_j,k}$ and some particular values of the zeta function
(Proposition \ref{z3}). Second, using the coefficients
$a_{\alpha_j,k}$, we
obtain a larger information on the values of the residues of the
associated zeta function (Proposition \ref{aaab3}).

\begin{prop} Let $S$ be a totally regular sequence of spectral type with
order $\alpha_N$ and genus $\ge$. Then (recall $\alpha_j\leq \ge$):
\begin{align*}
a_{\alpha_j,0}&=\left\{\begin{array}{ll}\Gamma(-\alpha_j)c_{\alpha_j,0}&
\alpha_j\not\in\Z~{\rm or}~\alpha_j\in \Z~{\rm and}~\alpha_j\leq
-1,\\
\frac{(-1)^{\alpha_j}}{\alpha_j}\hspace{-.05in}\left(\frac{1}{\alpha_j!}c_{\alpha_j,0}-\Rs{0}{\alpha_j}\zeta(s,S)\right)&\alpha_j\in
\Z~{\rm and}~0<\alpha_j\leq \ge,\\
-\Rs{0}{0}\zeta'(s,S)&\alpha_j=0,
\end{array}
\right.\\
a_{\alpha_j,1}&=\left\{\begin{array}{ll}0& \alpha_j\not\in\Z~{\rm
or}~\alpha_j\in \Z~{\rm and}~\alpha_j\leq -1,\\
\frac{(-1)^{\alpha_j+1}}{\alpha_j!}c_{\alpha_j,0}&\alpha_j\in
\Z~{\rm and}~0<\alpha_j\leq \ge,\\
-\Rs{0}{0}\zeta(s,S)=-c_{0,0}&\alpha_j=0.
\end{array}
\right.
\end{align*}
\label{z3}
\end{prop}

\begin{proof} We adapt here a technique introduced in \cite{Spr2}.
Integrate $\ge$ times with respect to $t$ in equation (\ref{f1}). We
obtain
\[
\zeta(s,S)=
\frac{(-1)^{\ge}s}{\Gamma(s-\ge)}\int_0^\infty t^{s-\ge-1}\frac{1}{2\pi
i}\int_{\Lambda_{\theta,c}}\frac{\e^{-\lambda
t}}{(-\lambda)^{\ge+1}} \log\Gamma(-\lambda,S) d\lambda dt,
\]

Next, split the $t$-integral at $t=1$. The $t$-integral $\int_1^\infty$ is a
regular function of $s$; in the $t$-integral $\int_0^1$, we change the
contour to $\Lambda_{\theta,-c}-C_c$, where $C_c$ is a circle with
center at the origin and radius $c$, and the $\lambda$ integration along $C_c$
vanishes, since by definition 
${{\rm Res}_{0}}_{_{_{\hspace{-.35in}\begin{array}{c}\scriptstyle{\lambda=0}\end{array}}\hspace{-.01in}}}
\frac{\e^{-\lambda t}}{(-\lambda)^\ge} \log
\Gamma(-\lambda,S)=0$ (in fact $\log\Gamma(-\lambda,S)$ has
a zero of order $\ge+1$ at $\lambda=0$).
Thus,
\[
\zeta(s,S) =\frac{(-1)^{\ge}s}{\Gamma(s-\ge)}\left(\int_0^1 t^{s-1}\frac{1}{2\pi
i}\int_{\Lambda_{\theta,-c}}\frac{\e^{-\lambda}}{(-\lambda)^{\ge+1}}
\log\Gamma(-\lambda/t,S) d\lambda dt +r(s)\right),
\]
where $r(s)$ is regular for all $s$ ($\Re(s)>\alpha_N-\epsilon$).
We use the expansion of the $\Gamma$-function and equations
(\ref{b1}) and (\ref{b2}) in the appendix, in order to obtain
\beq\label{ee1x}\begin{aligned}
\zeta(s,S)= \frac{(-1)^{\ge+1}s}{\Gamma(s-\ge)}&\left(
\sum_{j=0}^N 
\frac{1}{\Gamma(\ge+1-\alpha_j)}\frac{a_{\alpha_j,0}}{s-\alpha_j} 
+ \sum_{j=0}^N 
\frac{\psi(\ge+1-\alpha_j)}{\Gamma(\ge+1-\alpha_j)}
\frac{a_{\alpha_j,1}}{s-\alpha_j}
\right.\\
&\left.
+\sum_{j=0}^N 
\frac{1}{\Gamma(\ge+1-\alpha_j)}\frac{a_{\alpha_j,1}}{(s-\alpha_j)^2}+ r(s)\right).
\end{aligned}
\eeq

By direct inspection of the above equation, since the zeta function
can not have double poles and $S$ is totally regular, it follows that $a_{\alpha_j,1}=0$,
whenever $\alpha_j$ is not an integer, or $\alpha_j$ is an
integer and $\alpha_j\leq -1$ or $\alpha_j\geq \ge+1$. Now, if
$\alpha_j$ is not an integer or it is an integer and $\alpha_j\leq
-1$ or $\alpha_j\geq \ge+1$, then the zeta function has a pole at
$s=\alpha_j$ and
\[
\Ru_{s=\alpha_j}\zeta(s,S)=\frac{(-1)^{\ge+1}\alpha_j}{\Gamma(\alpha_j-\ge)}\frac{a_{\alpha_j,0}}{\Gamma(\ge+1-\alpha_j)}
=-\alpha_j \frac{\sin\alpha_j\pi}{\pi} a_{\alpha_j,0},
\]
therefore, since by Propositions \ref{coeff} and \ref{zeta1}
\[
\Ru_{s=\alpha_j}\zeta(s,S)=\frac{c_{\alpha_j,0}}{\Gamma(\alpha_j)}=-\alpha_j \frac{\sin\alpha_j\pi}{\pi} a_{\alpha_j,0},
\]
we obtain that $a_{\alpha_j,0}=\Gamma(-\alpha_j)c_{\alpha_j,0}$.
If $\alpha_j$ is an integer and $0\leq \alpha_j\leq \ge$, using the
expansion for $s$ near $\alpha_j=k$ of $\frac{s(s-1)\dots
(s-\ge)}{\Gamma(s)} =\frac{s}{\Gamma(s-\ge)}$ in equation (\ref{ee1x}), we obtain after some calculations
\[
\zeta(s,S)=(-1)^{k+1}(s-k+k)\left(\frac{a_{k,1}}{s-k}+a_{k,0}+O(s-k)\right).
\]

We distinguish two cases. First, if $k=0$, then near $s=0$
\[
\zeta(s,S)=-a_{0,1}-a_{0,0}s+O(s^2),
\]
as in Theorem \ref{zeta2}. Second, if $0<k\leq \ge$, then near $s=k$,
\[
\zeta(s,S)=(-1)^{k+1}k\frac{a_{k,1}}{s-k}+(-1)^{k+1}k a_{k,0}
+(-1)^{k+1}a_{k,1}+O(s-k),
\]
and using Proposition \ref{zeta1}, this completes the proof.
\end{proof}

\begin{prop} \label{aaab3} Let $S$ be a totally regular sequence of spectral type with order
$\alpha_N$ and genus $\ge$. Then, the associated zeta function
is regular in the complex half-plane $\Re(s)>\alpha_N-\epsilon$ up
to a finite set of at most $N+1$ simple poles located at
$s=\alpha_N,\alpha_{N-1},\dots, \alpha_0$.

\begin{align*}
\Ru_{s=\alpha_j}\zeta(s,S)&=\frac{a_{\alpha_j,0}}{\Gamma(\alpha_j)\Gamma(-\alpha_j)},
\hspace{.4in} \alpha_j\not=-[|\alpha_N|],\dots,\ge,\\
\Ru_{s=\alpha_j}\zeta(s,S)&=\left\{\begin{array}{ll}(-1)^{\alpha_j+1}\alpha_j
a_{\alpha_j,1},& \alpha_j=1,\dots,\ge,\\
0,&\alpha_j=-[|\alpha_N|],\dots,0,\end{array}\right.\\
\Rz_{s=\alpha_j}\zeta(s,S)&=\left\{\begin{array}{ll}(-1)^{\alpha_j+1}\alpha_j
a_{\alpha_j,0},&\alpha_j=-[|\alpha_N|],\dots,-1,\\
(-1)^{\alpha_j+1}\left(\alpha_j
a_{\alpha_j,0}+a_{\alpha_j,1}\right),&\alpha_j=1,\dots,\ge,\\
- a_{\alpha_j,1}&\alpha_j=0.\end{array}\right.
\end{align*}
\end{prop}

We conclude this section with some properties of sums and
multiplication by a scalar of sequences of spectral type. All
the proofs are straightforward, up to \ref{lin4}, where we can use
the Young inequality.
Let  $S=\{\lambda_{n}\}_{n\in \N_0}$ be a sequence and $y$  a
positive  number. Denote by $yS$ the sequence
$yS=\{y\lambda_{n}\}_{n\in \N_0}$. Let $n=(n_1,\dots, n_I)$ be an integer vector. If
$S_i=\{\lambda_{(i),n_i}\}_{n_i\in \N_0}$, $i=1,\dots, I$, is  a
finite set of sequences, we use the notation $\sum_{i=1}^I S_{(i)}$
for the sum of the sequences, namely the sequence
$\left\{\lambda_n=\sum_{i=1}^I
\lambda_{(i),n_1}\right\}_{n\in (\N_0)^I}$.

\begin{lem}\label{lin1}
For any positive real $y$, $S=\{\lambda_{n}\}_{n\in \N_0}$ is  a
sequence with finite exponent $\ec$, genus $\ge$ and spectral sector
$\Sigma_{\theta,c}$,  if and only if $yS$ is  a  sequence with
finite exponent $\ec$, genus $\ge$, and spectral sector
$\Sigma_{\theta,yc}$. $S$ is of spectral type of order $\alpha_N$ if
and only if $yS$  is of spectral type of order $\alpha_N$; $S$ is
regular, simply or totally regular if and only if $yS$ is regular,
simply or totally regular respectively.
\end{lem}

\begin{lem}\label{lin3}
If $S_{(i)}=\{\lambda_{(i),n_i}\}_{n_i\in \N_0}$, $i=1,\dots, I$, is
a finite family of sequences, then
\[
f\left(t,\sum_{i=1}^I
S_{(i)}\right)=1+\prod_{i=1}^I(f(t,S_{(i)})-1).
\]
\end{lem}

\begin{lem}\label{lin4}
If $S_{(i)}=\{\lambda_{(i),n_i}\}_{n_i\in \N_0}$, $i=1,\dots, I$,
is a finite family of sequences of spectral type of finite exponents
$\ec_{(i)}$, genus $\ge_{(i)}$, asymptotic domains
$AD_{\theta_{(i)},c_{(i)}}$, and orders $\alpha_{(i),N_{(i)}}$,
then the sum sequence $S=\sum_{i=1}^I S_{(i)}$ is a sequence of
spectral type with exponent $\ec_{(0)}=\sum_{i=1}^I \ec_{(i)}$,
asymptotic domain $AD_{\theta,c}$, where $\theta={\rm
max}(\theta_{(i)})$ and $c={\rm min}(c_{(i)})$, and order
$\alpha_N={\rm max}(\alpha_{(i),N{(i)}})$.
\end{lem}
\begin{lem}\label{lin5}
If $S_{(i)}=\{\lambda_{(i),n_i}\}_{n_i\in \N_0}$, $i=1,\dots, I$,
is a finite family of simply/totally regular sequences of spectral
type, then the sum sequence $S=\sum_{i=1}^I S_{(i)}$ is a
simply/totally regular sequence of spectral type.
\end{lem}

\section{Spectral decomposition}
\label{s2}

The aim of this section is to obtain the main zeta
invariants of double sequences.  It is clear that using multi indices the results extend to multiple sequences. Our strategy to deal with double
sequences is to define a class of double sequences that can be
decomposed as sums of simple sequences relatively to some fixed
simple sequence, in a suitable way. The double sequences of this
class are said to be spectrally decomposable, and are defined in the
following Definition \ref{spdec}. The idea of this type of
decomposition was suggested by the results of Br\"uning and Seeley (in
particular see \cite{BS2}), and has already been applied with
success in a number of cases \cite{Spr3}  \cite{Spr5} \cite{Spr7}. For
the zeta function associated to spectrally decomposable double
sequences, we obtain conditions for regularity at $s=0$ and a
formula that gives the values of $\zeta(0)$ and $\zeta'(0)$ (Theorem
\ref{sdl}).

Let $S=\{\lambda_{n,k}\}_{n,k=1}^\infty$ be a double sequence of non
vanishing complex numbers with unique accumulation point the
infinity, finite exponent $s_0=\ec(S)$ and genus $p_0=\ge(S)$. Assume if necessary that the element of $S$ are ordered as $0<|\lambda_{1,1}|\leq|\lambda_{1,2}|\leq |\lambda_{2,1}|\leq \dots$. We use the notation $S_n$ ($S_k$) to denote the simple sequence with fixed $n$ ($k$). We call the exponents of $S_n$ and $S_k$ the relative exponents of $S$, and we use the notation $(s_0=\ec(S),s_1=\ec(S_k),s_2=\ec(S_n))$. We define relative genus accordingly.

\begin{defi} Let $S\hspace{-1pt}=\hspace{-1pt}\{\lambda_{n,k}\}_{n,k=1}^\infty$ be a double
sequence with finite exponents $\hspace{-.5pt}(\hspace{-.5pt}s_0,s_1,s_2)$, genus
$(p_0,p_1,p_2)$, and positive spectral sector
$\Sigma_{\theta_0,c_0}$. Let $U=\{u_n\}_{n=1}^\infty$ be a totally
regular sequence of spectral type of order $u'\leq 0$ with exponent
$r_0$, genus $q$, domain $AD_{\phi,d}$. We say that $S$ is
{\it spectrally decomposable over $U$ with power $\kappa$, length $\ell$ and
asymptotic domain $AD_{\theta,c}$, with $c={\rm min}(c_0,d,c')$,
$\theta={\rm max}(\theta_0,\phi,\theta')$}, if there exist positive
real numbers $\kappa$, $\ell$ (integer), $c'$, and $\theta'$, with
$0< \theta'<\pi$,   such that:
\begin{enumerate}
\item the sequence
$\widetilde{S}_n=u_n^{-\kappa}S_n
=\left\{\frac{\lambda_{n,k}}{u^\kappa_n}\right\}_{k=1}^\infty$ has
spectral sector $\Sigma_{\theta',c'}$, and is a totally regular
sequence of spectral type of order $\leq 0$ for each
$n$;
\item the logarithmic $\Gamma$-function associated to  $\widetilde{S}_n$ has an asymptotic expansion  for large
$n$ uniformly in $\lambda$ for $\lambda$ in
$AD_{\theta,c}$, of the following form
\[
\log\Gamma(-\lambda,u_n^{-\kappa} S_n)=\sum_{h=0}^{\ell}
\phi_{\sigma_h}(\lambda) u_n^{-\sigma_h}+\sum_{l=0}^{L}
P_{\rho_l}(\lambda) u_n^{-\rho_l}\log u_n+o(u_n^{-r_0}),
\]
where $\sigma_h$ and $\rho_l$ are real numbers with $\sigma_0<\dots <\sigma_\ell$, $\rho_0<\dots <\rho_L$, the
$P_{\rho_l}(\lambda)$ are polynomials in $\lambda$ satisfying the
condition $P_{\rho_l}(0)=0$, $\ell$ and $L$ are the larger integers
such that $\sigma_\ell\leq r_0$ and $\rho_L\leq r_0$,  and $\sigma_0$ and $\rho_0$ satisfy the
condition $u'<{\rm min}(\sigma_0,\rho_0-1)$.
\end{enumerate}
\label{spdec}
\end{defi}

We point out that the condition $P_{\rho_l}(0)=0$ is introduced in Definition \ref{spdec} expressely in order to avoid the dependence of the final result on the logarithmic terms. More precisely, without this condition, the function $A_{0,0}(s)$, introduced in Lemma \ref{lp2}, would depend on the values of the coefficients of the polynomials $P_{\rho_l}(\lambda)$. As a consequence, also the coefficients of the expansion of the zeta function at zero given in Theorem \ref{sdl} would depend on the coefficients of the $P_{\rho_l}(\lambda)$.

The proof of the following lemma is  essentially based on
standard properties of double limits.

\begin{lem} \label{pippo2} The functions $\phi_{\sigma_h}(\lambda)$
appearing in Definition \ref{spdec} have an asymptotic expansion of
order $\leq 0$ for large $\lambda$ in $AD_{\theta,c}$ with respect
to the asymptotic sequence $\{(-\lambda)^\alpha \log^k
(-\lambda)\}_\alpha$, $k=0,1$, where the unique logarithmic terms are
of the form $(-\lambda)^m\log (-\lambda)$, with integer $m$ such
that $m\leq p_2$. The maximum order of the
polynomials $P_{\rho_l}(\lambda)$ is lower or equal to $p_2$.
\end{lem}

When $S$ is a double sequence, we consider the general case where we
do not know an explicit expansion for $\log\Gamma(-\lambda,S)$, and
hence we can not apply Theorem \ref{zeta2}. On the other side we
suppose that we do know such expansions for the sequences
$\widetilde{S}_n$, and our aim is to use this information to obtain
the result for $S$. We could write
$\log\Gamma(-\lambda,S)=\sum_{n=1}^\infty \log\Gamma(-\lambda,S_n)$,
and use this decomposition in the formula (\ref{f1}) for
$\zeta(s,S)$. Then, we could proceed as in the proof of Proposition
\ref{z3} using the known
information on the expansion of $S_n$, and eventually 
try to sum up on $n$. This procedure does not work. The series in
$n$ may not converge. To overcome this difficulty, the idea is to
insert an $s$ in the general term. In other words, we decompose the zeta function as: 
$\zeta(s,S)=\sum_{n=1}^\infty u_n^{-\kappa s}
\zeta(s,u_n^{-\kappa}S_n)$ (see Lemma \ref{lp1} below). The  power $\kappa$ is necessary to guarantee the
existence of an asymptotic expansion of
$\log\Gamma(-\lambda,\widetilde{S}_n)$ for large $n$.

\begin{lem} Let $S$ be spectrally decomposable over $U$ as in Definition \ref{spdec}. Then, there
exist constants $N$, $\alpha_N<\dots<\alpha_0$, $a_{\alpha_j,k,n}$ and
$b_{\sigma_h,\alpha_j,k}$, with $\alpha_N\leq 0$, such that
\[
\log\Gamma(-\lambda,\widetilde{S}_n)=\sum_{j=0}^N\sum_{k=0}^1
a_{\alpha_j,k,n}
(-\lambda)^{\alpha_j}\log^k(-\lambda)+o((-\lambda)^{\alpha_N}),
\]
for all $n$, and
\[
\phi_{\sigma_h}(\lambda)=\sum_{j=0}^N \sum_{k=0}^1
b_{\sigma_h,\alpha_j,k} (-\lambda)^{\alpha_j}\log^k
(-\lambda)+o((-\lambda)^{\alpha_N}),
\]
for all $\sigma_h$, for large $\lambda$ uniformly in
$AD_{\theta,c}$, and where
$a_{\alpha_j,1,n}=b_{\sigma_h,\alpha_j,1}=0$ for all
$\alpha_j\not= 0,1, \dots, p_2$.
\label{lll1}
\end{lem}
\begin{proof} Since $\widetilde{S}_n=u_n^{-\kappa}S_n$ is of spectral type
for each $n$, we have that the sequences of powers $\alpha_j$
appearing in the asymptotic expansion of the $\Gamma$-functions are
all upper bounded by $s_1$ (by Proposition \ref{zeta1}). Hence, they
can all only accumulate at $-\infty$, and hence we can collect them
in a unique sequence $\{\alpha_j\}_{j=0}^N$, starting at the bigger
non vanishing one and where $\alpha_N$ is the smaller of the order
of the $\widetilde{S}_n$. Notice also that since the expansions of
the $\log\Gamma(-\lambda,S_n)$ are all of order $\leq 0$, it follows
that $\alpha_N\leq 0$. For the domain, if the condition in
Definition \ref{spdec} are satisfied, all the $\widetilde{S}_n$ are
contained in the positive spectral sector $\Sigma_{\theta,c}$. The
same argument works for the $\phi_{\sigma_h}$, and again we can
reset $\alpha_N$ to be the smaller of the orders and it will always
be negative or null. The condition on the coefficients follows from
the point (1) (the fact that $\widetilde{S}_n$ is totally regular)
of Definition \ref{spdec} (see also Remark \ref{form}). \end{proof}

\begin{lem}
\label{lp1} Let $S$ be spectrally decomposable over $U$ as in Definition \ref{spdec}. Then,
we have the following contour integral representation for the associated
zeta function 
\[
\zeta(s,S)=\frac{s}{\Gamma(s)}\int_0^\infty t^{s-1}\frac{1}{2\pi
i}\int_{\Lambda_{\theta,c}}\frac{\e^{-\lambda t}}{-\lambda}
\T(s,\lambda,S,U) d\lambda,
\]
where
\[
\T(s,\lambda,S,U)= \sum_{n=1}^\infty u_n^{-\kappa s}
\log\Gamma(-\lambda,u_n^{-\kappa} S_n).
\]
\end{lem}
\begin{proof} By uniform convergence of the series, for $\Re(s)>s_0$, we can write
\[
\zeta(s,S)=\sum_{n=1}^\infty u_n^{-\kappa s}
\zeta(s,u_n^{-\kappa}S_n),
\]
where 
\[
\zeta(s,u_n^{-\kappa s}S_n)=u_n^{-\kappa s}\sum_{k=1}^\infty \lambda_{n,k}^{-s},
\] 
is well defined by the series since
$s_0\geq s_2$. Now, applying equation (\ref{f1b}) of Lemma \ref{f} to $\zeta(s,u_n^{-\kappa s}S_n)$, 
when $\Re(s)>s_0\geq
s_1,s_2$, thanks to uniform convergence of the integral, we obtain the thesis.
\end{proof}

\begin{lem} \label{lp2} Let $S$ be spectrally decomposable over $U$ as in Definition \ref{spdec}.
Then, the function $\T(s,\lambda,S)$ can be extended analytically to
the half plane $\Re(s)> -\epsilon$ (positive small $\epsilon$) by
the following formula for all $\lambda\in AD_{\theta,c}$,
\[
\T(s,\lambda,S,U)=\P(s,\lambda,S,U)+\sum_{h=0}^\ell
\phi_{\sigma_h}(\lambda)\zeta(\kappa s+\sigma_h,U)-\sum_{l=0}^L
P_{\rho_l}(\lambda)\zeta'(\kappa s+\rho_l,U),
\]

Moreover, the function $\P(s,\lambda,S,U)$ is regular in $s$ for $\Re(s)>-\epsilon$ and has the following
expansion for large $\lambda$ uniformly in $AD_{\theta,c}$
\[
\P(s,\lambda,S,U)=\sum_{j=0}^N A_{\alpha_j,0}(s)(-\lambda)^{\alpha_j}+
\sum_{j=0}^{p_2}
A_{j,1}(s)(-\lambda)^{j}\log(-\lambda)+o((-\lambda)^{\alpha_N}),
\]
where $N$ is defined in Lemma \ref{lll1} and the coefficients
\begin{align*}
A_{\alpha_j,0}(s)&=\sum_{n=1}^\infty \left(a_{\alpha_j, 0,n}
-\sum_{h=0}^\ell b_{\sigma_h,\alpha_j,0}u_n^{-\sigma_h}
\right)u_n^{-\kappa s}, ~~~~\alpha_j\not=0,1,\dots,p_2,\\
A_{0,0}(s)&=\sum_{n=1}^\infty \left(a_{0, 0,n} -\sum_{h=0}^\ell
b_{\sigma_h,0,0}u_n^{-\sigma_h}\right)u_n^{-\kappa s},\\
A_{j,0}(s)&=\sum_{n=1}^\infty \left(a_{j, 0,n} -\sum_{h=0}^\ell
b_{\sigma_h,j,0}u_n^{-\sigma_h}\hspace{-2pt} -\hspace{-2pt}\sum_{l=0}^L
p_{\rho_l,j}u_n^{-\rho_l}\log u_n \hspace{-3pt}\right)\hspace{-2pt} u_n^{-\kappa s},~1\leq j\leq p_2,\\
A_{j,1}(s)&=\sum_{n=1}^\infty \left(a_{j, 1,n} -\sum_{h=0}^\ell
b_{\sigma_h,j,1}u_n^{-\sigma_h}\right)u_n^{-\kappa s},
~~~0\leq j\leq p_2,
\end{align*}
are regular functions of $s$ for $\Re(s)>-\epsilon$.
\end{lem}

Note that $\T$
and $\P$ are regular does not mean that $\zeta(s,S)$ is regular,
just that we need not to bother about the dependence on $s$ coming
from $\T$.

\begin{proof} We would like to expand $\T$ for large $\lambda$ in order to
proceed as in the proof of Proposition \ref{z3}. We could use the
expansion of $\log\Gamma(-\lambda,u_n^{-\kappa} S_n)$ for large
$\lambda$ to get that of $\T$. Unfortunately, this does not work for
the following reason. When we expand $\T$ and we perform the
integrals in $\lambda$ and $t$ for each term, the resulting
functions have poles in $s$ at $s=0$ and the sum over $n$ also gives
a pole in the same point, thus we get some terms with a double pole,
and hence we can not use the formula (\ref{f1}) to compute the
derivative in $s$ at $s=0$. What we can do, is to split in two
terms, each one having only simple poles. We do this as follows.
Since $S$ is spectrally decomposable over $U$ with power $\kappa$,
we have the expansion (observe that
$\log\Gamma(-\lambda,\widetilde{S}_n)=\log\Gamma(-\lambda
u_n^\kappa,S_n)$)
\[
\log\Gamma(-\lambda, \widetilde{S}_n)=\sum_{h=0}^\ell
\phi_{\sigma_h}(\lambda) u_n^{-\sigma_h}+\sum_{l=0}^L
P_{\rho_l}(\lambda) u_n^{-\rho_l}\log u_n+o(u_n^{-K}).
\]

This allows to split $\T$ in two terms
\begin{align*}
\T(s,\lambda,S,U)&=\P(s,\lambda,S,U)+\sum_{h=0}^\ell \sum_{n=1}^\infty
\phi_{\sigma_h}(\lambda) u_n^{-\kappa s-\sigma_h}+\sum_{l=0}^L
P_{\rho_l}(\lambda) u_n^{-\kappa s-\rho_l}\log u_n\\
&= \P(s,\lambda,S,U)+\sum_{h=0}^\ell
\phi_{\sigma_h}(\lambda)\zeta(\kappa s+\sigma_h,U)-\sum_{l=0}^L
P_{\rho_l}(\lambda) \zeta'(\kappa s+\rho_l,U),
\end{align*}
where
\[
\P(s,\lambda,S,U)=\sum_{n=1}^\infty  u_n^{-\kappa s}\hspace{-2pt}\left(\hspace{-2pt}\log\Gamma(-\lambda
u_n^{\kappa}, S_n)\hspace{-2pt}-\hspace{-2pt}\sum_{h=0}^\ell \phi_{\sigma_h}(\lambda)
u_n^{-\sigma_h}\hspace{-2pt}-\hspace{-2pt}\sum_{l=0}^L
P_{\rho_l}(\lambda)u_n^{-\rho_l}\log u_n\hspace{-2pt} \right)\hspace{-4pt}.
\]

Since, by Definition \ref{spdec}, the term in brackets in the above sum is $O(u_n^{-r_0-\epsilon})$ uniformly in
$\lambda$ for $\lambda\in AD_{\theta,c}$, it follows that
\[
\P(s,\lambda,S,U)<\sum_{n=1}^\infty u^{-\kappa s-r_0-\epsilon}_n,
\]
and hence $\P$ is regular when $\Re(\kappa s+r_0+\epsilon)\geq r_0$, i.e.
$\Re(s)\geq -\frac{\epsilon}{\kappa}$. As $\kappa>0$, $\P$ is regular at $s=0$, uniformly for
$\lambda$ in $AD_{\theta,c}$. Thus we can get the expansion of $\P$
as follows. It is clear that the expansion for large $\lambda$ of
$\log\Gamma(-\lambda, S_n)$ gives that of $\log\Gamma(-\lambda
u_n^{\kappa}, S_n)$ ($n$ fixed). Thus we can use the expansions
given in Lemma \ref{lll1}. Substituting in the definition of $\P$ we
have
\begin{align*}
\P(s,\lambda,S,U)=\sum_{n=1}^\infty u_n^{-\kappa s}\sum_{j=0}^N
\sum_{k=0}^1 &\left(a_{\alpha_j,k,n} (-\lambda)^{\alpha_j}\log^k
(-\lambda)-\sum_{l=0}^L P_{\rho_l}(\lambda)u_n^{-\rho_l}\log u_n\right.\\
&\left.
-\sum_{h=0}^\ell b_{\sigma_h,\alpha_j,k}u_n^{-\sigma_h}(-\lambda)^{\alpha_j}\log^k(-\lambda)
\right)\hspace{-3.1pt}+o((-\lambda)^{\alpha_N}),
\end{align*}
since the series defining $\P$ is uniformly convergent for
$\Re(s)>s_2$,
\begin{align*}
\P(s,\lambda,S,U)=&\sum_{n=1}^\infty u_n^{-\kappa s}\sum_{j=0}^N \left(
a_{\alpha_j,0,n} (-\lambda)^{\alpha_j} -\sum_{h=0}^\ell
b_{\sigma_h,\alpha_j,0}u_n^{-\sigma_h}(-\lambda)^{\alpha_j} \right)\\
&+\sum_{n=1}^\infty u_n^{-\kappa s}\sum_{j=0}^{p_2} \left( a_{j,1,n}
(-\lambda)^j\log(-\lambda)
-\sum_{l=0}^L P_{\rho_l}(\lambda)u_n^{-\rho_l}\log u_n \right.\\
&\hspace{72pt}\left.-\sum_{h=0}^\ell
b_{\sigma_h,j,1}u_n^{-\sigma_h}(-\lambda)^j\log (-\lambda)\right)+
o((-\lambda)^{\alpha_N}),
\end{align*}
writing (recall $P_{\rho_l}(0)=0$)
$P_{\rho_l}(\lambda)=\sum_{m=1}^{G_{\rho_l}} p_{\rho_l,m}
(-\lambda)^m$, we have
\[
\sum_{l=0}^L P_{\rho_l}(\lambda) \zeta'(\kappa s+\rho_l,U)=
\sum_{m=1}^G\left(\sum_{l=0}^L p_{\rho_l,m} \zeta'(\kappa
s+\rho_l,U)\right)(-\lambda)^m,
\]
where $G={\rm max}(G_{\rho_l})\leq p_2$, and we set $p_{\rho_l,m}=0$
for $m>G_{\rho_l}$; then
\begin{align*}
\P(s,\lambda,S,U)=&\sum_{n=1}^\infty \left(a_{\alpha_j,0,n}
-\sum_{h=0}^\ell b_{\sigma_h,\alpha_j,0}u_n^{-\sigma_h}
\right)u_n^{-\kappa s}(-\lambda)^{\alpha_j}\\
&-\sum_{n=1}^\infty\sum_{m=1}^G\left(\sum_{l=0}^L p_{\rho_l,m}
u_n^{-\kappa s-\rho_l}\log u_n\right)(-\lambda)^m\\
&+\hspace{-2pt}\sum_{n=1}^\infty \sum_{j=0}^{p_2}\hspace{-2pt}\left(\hspace{-2pt}a_{j,1,n}
\hspace{-2pt}-\hspace{-2pt}\sum_{h=0}^\ell b_{\sigma_h,j,1}u_n^{-\sigma_h}\hspace{-2pt}\right)\hspace{-2pt}u_n^{-\kappa
s}(-\lambda)^j\log(-\lambda)+o((-\lambda)^{\alpha_N}).
\end{align*}

This gives the coefficients stated in the thesis and since $\P(s,\lambda,S,U)$
is regular at $s=0$, this also shows that they are regular at $s=0$. \end{proof}

\begin{lem} Let $S$ be spectrally decomposable over $U$ as in Definition \ref{spdec}. Then, near
$s=0$, we have the following representation for the zeta function
\[
\Gamma(s)\zeta(s,S)
=\gamma
A_{0,1}(s)-A_{0,0}(s)-\frac{A_{0,1}(s)}{s}
+s\sum_{h=0}^{\ell}\Phi_{\sigma_h}(s)\zeta(\kappa s+\sigma_h,U)+s
r(s),
\]
where $r(s)$ is regular near $s=0$ and 
\[
\Phi_{\sigma_h}(s)=\int_0^\infty t^{s-1}\frac{1}{2\pi
i}\int_{\Lambda_{\theta,c}}\frac{\e^{-\lambda t}}{-\lambda}
\phi_{\sigma_h}(\lambda)d\lambda d t.
\]
\label{lp3}
\end{lem}
\begin{proof} Let $\Re(s)>s_0$. By Lemma \ref{lp1}
\[
\zeta(s,S)=\frac{s}{\Gamma(s)}\int_0^\infty t^{s-1}\frac{1}{2\pi
i}\int_{\Lambda_{\theta,c}}\frac{\e^{-\lambda t}}{-\lambda}
\T(s,\lambda,S) d\lambda.
\]
Substituting the decomposition of $\T$ given by Lemma \ref{lp2}, we obtain
\[
\zeta(s,S)=
\zeta_1(s)+\zeta_2(s)+\zeta_3(s),
\]
\begin{align*}
\zeta_1(s)&=\frac{s}{\Gamma(s)}\int_0^\infty
t^{s-1}\frac{1}{2\pi i}\int_{\Lambda_{\theta,c}}\frac{\e^{-\lambda
t}}{-\lambda} \P(s,\lambda,S,U) d\lambda d t\\
\zeta_2(s)&=\frac{s}{\Gamma(s)} \sum_{h=0}^H \zeta(\kappa s+\sigma_h,U)
\int_0^\infty t^{s-1}\frac{1}{2\pi
i}\int_{\Lambda_{\theta,c}}\frac{\e^{-\lambda t}}{-\lambda}
\phi_{\sigma_h}(\lambda) d\lambda d t\\
\zeta_3(s)&=-\frac{s}{\Gamma(s)} \sum_{l=0}^{L} P_{\rho_l}(\lambda)
\zeta'(\kappa s +\rho_l,U)\int_0^\infty t^{s-1}\frac{1}{2\pi
i}\int_{\Lambda_{\theta,c}}\frac{\e^{-\lambda t}}{-\lambda}
P_{\rho_l}(\lambda) d\lambda d t.
\end{align*}

Now, $\zeta_2(s)$ can
immediately be written as stated in the thesis, while $\zeta_3(s)$ vanishes by direct calculation (using the integral given in the appendix).  
For $\zeta_1(s)$, we can use
the expansion given in Lemma \ref{lp2} for $\P$ as follows. First,
split the integral ($\Re(s)>s_0$)
\begin{align*}
\zeta_1(s)=&\frac{s}{\Gamma(s)}\int_0^1 t^{s-1}\frac{1}{2\pi
i}\int_{\Lambda_{\theta,c}}\frac{\e^{-\lambda t}}{-\lambda}
\P(s,\lambda,S,U) d\lambda d t\\
&+\frac{s}{\Gamma(s)}\int_1^\infty t^{s-1}\frac{1}{2\pi
i}\int_{\Lambda_{\theta,c}}\frac{\e^{-\lambda t}}{-\lambda}
\P(s,\lambda,S,U) d\lambda d t.
\end{align*}

For the second term, note that $\P$ is regular for $\Re(s)\geq 0$,
thus we can extend the $s$-domain to the non negative $s$-plane;
also $\P$ diverges at most like a power in $AD_{\theta,c}$, and
hence the $\lambda$ integral is bounded due to the presence of the
$\e^{-\lambda t}$ ($t>0$); eventually, the $t$ integral is also
bounded by the same factor for all $s$. This means that the second
term is the product of $\frac{s}{\Gamma(s)}$ with a function
$r_1(s)$ that is regular for all $\Re(s)\geq 0$. For the first term,
we need to rewrite the contour as
$\Lambda_{\theta,c}=\Lambda_{\theta,-c}-C_c$. In the integral on the
new contour $\Lambda_{\theta,-c}$ we can use the expansion of $\P$,
the other gives
\[
\frac{1}{2\pi i}\int_{C_c}\frac{\e^{-\lambda t}}{-\lambda}
\P(s,\lambda,S,U) d\lambda=\P(s,0,S,U)=0,
\]
since $\T(s,\lambda,S,U)=0$ because $\log\Gamma(0,S_n)=0$ by definition, $P_{\rho_l}(0)=0$, and $\phi_{\sigma_h}(0)=0$ since
$\lambda=0$ belongs to $AD_{\theta,c}$.
We are left with
\begin{align*}
\zeta_1(s)&=\frac{s}{\Gamma(s)}\int_0^1 t^{s-1}\frac{1}{2\pi
i}\int_{\Lambda_{\theta,-c}}\frac{\e^{-\lambda t}}{-\lambda}
\P(s,\lambda,S,U) d\lambda d t,
\end{align*}
where we can use the expansion of $\P$, given in Lemma \ref{lp2}. We obtain (use equations (\ref{b1}) and (\ref{b2}) in the appendix)
\begin{align*}
\zeta_1(s)=&-\frac{s}{\Gamma(s)}\sum_{j=0}^N
A_{\alpha_j,0}(s)\frac{1}{s-\alpha_j}\frac{1}{\Gamma(1-\alpha_j)}
-\frac{s}{\Gamma(s)}\sum_{j=0}^1
A_{j,1}(s)\frac{1}{s-j}\frac{\psi(1-j)}{\Gamma(1-j)}\\
&- \frac{s}{\Gamma(s)}\sum_{j=0}^1
A_{j,1}(s)\frac{1}{(s-j)^2}\frac{1}{\Gamma(1-j)}+
\frac{s}{\Gamma(s)} r_2(s),
\end{align*}
with $r_2(s)$ regular near $s=0$. Since by definition and Lemma \ref{lll1}, $\alpha_N\leq 0$, near
$s=0$ the only non vanishing terms are those with $\alpha_j=0$,
hence
\[
\zeta_1(s)=\frac{s}{\Gamma(s)}\left(
-\frac{1}{s}A_{0,0}(s)+\frac{1}{s}\gamma
A_{0,1}(s)-\frac{1}{s^2}A_{0,1}(s) +r_3(s)\right),
\]
where $r_3(s)$ is regular near $s=0$. Summing up we have the thesis.
\end{proof}

\begin{lem}\label{ultimo} The functions $\Phi_{\sigma_h}(s)$ have
at most a double pole at $s=0$.
\end{lem}
\begin{proof} Just proceed as in the proof of Lemma \ref{lp3} and use the
asymptotic expansion given for $\phi_{\sigma_h}$ in Lemma
\ref{lll1}. \end{proof}

\begin{rem} \label{rp1} Note that $\zeta_2(s)$ can have at most a
simple pole at $s=0$. In fact, by definition $\zeta(s,U)$ has at
most simple poles. On the other side, $\Phi_{\sigma_h}(s)$ has
also at most double poles, by the previous lemma, and this implies
the statement, since $\frac{s}{\Gamma(s)}$ as a zero of order 2 at
$s=0$. Note also that we can not use the expansion of
$\phi_{\sigma_h}(\lambda)$ to compute the derivative of
$\zeta_2(s)$ at $s=0$, since the remainder of such an expansion
would appear in $\zeta_2'(0)$. We need explicit computation of the
integral defining the function $\Phi_{\sigma_h}(s)$ and an
expansion of $\Phi_{\sigma_h}(s)$ near $s=0$ at least up to the
constant term.
\end{rem}

\begin{teo} {\rm (Spectral decomposition lemma)}\label{sdl} Let the double sequence $S$ be spectrally decomposable over $U$ with power $\kappa$ and length $\ell$, then
\begin{align*}
\Ru_{s=0}\zeta(s,S)=&\frac{1}{\kappa}\sum_{h=0}^\ell
\Rd_{s=0}\Phi_{\sigma_h}(s) \Ru_{s=\sigma_h}\zeta(s,U),\\
\Rz_{s=0}\zeta(s,S)=&\sum_{h=0}^\ell\Rd_{s=0}\Phi_{\sigma_h}(s)\Rz_{s=\sigma_h}\zeta(s,U)-A_{0,1}(0)\\
&+\frac{1}{\kappa}\sum_{h=0}^\ell\Ru_{s=\sigma_h}\zeta(s,U)\left(\Ru_{s=0}\Phi_{\sigma_h}(s)+
\gamma\Rd_{s=0}\Phi_{\sigma_h}(s)\right),\\
\Rz_{s=0}\zeta'(s,S)
=&\frac{1}{\kappa}\left(\frac{\gamma^2}{2}-\frac{\pi^2}{12}\right)\sum_{h=0}^\ell
\Rd_{s=0}\Phi_{\sigma_h}(s)\Ru_{s=\sigma_h}\zeta(s,U)\\
&\hspace{-2.8pt}+\hspace{-2pt}\frac{\gamma}{\kappa}\sum_{h=0}^\ell
\Ru_{s=0}\Phi_{\sigma_h}(s)\Ru_{s=\sigma_h}\zeta(s,U)\hspace{-2pt}
+\gamma\hspace{-2pt}\sum_{h=0}^\ell
\Rd_{s=0}\Phi_{\sigma_h}(s)\Rz_{s=\sigma_h}\zeta(s,U)\\
&\hspace{-2.8pt}+\hspace{-2pt}\frac{1}{\kappa}\sum_{h=0}^\ell
\Rz_{s=0}\Phi_{\sigma_h}(s)\Ru_{s=\sigma_h}\zeta(s,U)\hspace{-2pt}
+\kappa\hspace{-2pt}\sum_{h=0}^\ell
\Rd_{s=0}\Phi_{\sigma_h}(s)\Rz_{s=\sigma_h}\zeta'(s,U)\\
&\hspace{-2.8pt}+\sum_{h=0}^\ell
\Ru_{s=0}\Phi_{\sigma_h}(s)\Rz_{s=\sigma_h}\zeta(s,U)
-A_{0,0}(0)-A'_{0,1}(0).
\end{align*}
\end{teo}
\begin{proof} Let $s$ be in a small neighborhood of $s=0$. Then, by Lemma \ref{lp3}
\beq\label{xxx}
\zeta(s,S)=\frac{1}{\Gamma(s)}\hspace{-1.4pt}\left(\hspace{-2pt}\gamma
A_{0,1}(s)-A_{0,0}(s)-\frac{A_{0,1}(s)}{s}
+s\sum_{h=0}^{\ell}\Phi_{\sigma_h}(s)\zeta(\kappa s+\sigma_h,U)+s
r(s)\hspace{-2pt}\right)\hspace{-2.6pt},
\eeq
where $r(s)$ is regular. We are able to expand $\zeta(s,S)$ near $s=0$, using the expansions of the single factors provided by the previous results. By Lemma \ref{lp2} the functions $A_{i,j}(s)$ are regular near $s=0$. By definition \ref{spdec}, the function $\zeta(s,U)$ has at most simple poles, and therefore we have the expansion
\[
\zeta(\kappa s+\sigma_h,U)=\frac{1}{\kappa}\Ru_{s=\sigma_h}\zeta(s,U)\frac{1}{s}+\Rz_{s=\sigma_h}\zeta(s,U)+\kappa\Rz_{s=\sigma_h}\zeta'(s,U)s+O(s^2),
\]
near $s=0$. By Lemma \ref{ultimo} and Remark \ref{rp1}, the functions $\Phi_{\sigma_h}(s)$ have at most a pole of order 2 at $s=0$, and therefore we have an expansion of the form
\[
\Phi_{\sigma_h}(s)=\frac{\Rd_{s=0}\Phi_{\sigma_h}(s)}{s^2}+\frac{\Ru_{s=0}\Phi_{\sigma_h}(s)}{s}+\Rz_{s=0}\Phi_{\sigma_h}(s)+O(s).
\]

Substituting these expansions and the classical expansion of $1/\Gamma(s)$ in equation (\ref{xxx}), after some calculations, we obtain the thesis.
\end{proof}

The result looks a bit hard as it stands, but this is a consequence of its large generality. The main point is precisely that this result is very general and applies in a very large number of cases. It turns out that in applications the explicit form of the different terms becomes nicer and many simplifications arise, but in order to cover the various situations the general setting is necessary. 

\section{Sums of sequences of spectral type}
\label{s3}

Let $S_{i}$, $i=1,\dots, I$ be a  finite set of totally regular
sequences of spectral type.  In this section we study the
decomposition properties of the sum of these sequences with respect to
one of them, and applying Theorem \ref{sdl} we obtain the expansion of the zeta function associated to the sum of the sequences. By the results at the end of Section \ref{s1} on
the sum of sequences, it is clear that  we can assume $I=2$ without loss of generality.
We will use the notation with a subscript index between parenthesis $X_{(i)}$ to denote the sequence to which the quantity $X$ is referred.

Let $S_{(i)}=\{\lambda_{(i),n_i}\}_{n_i\in \N_0}$, $i=1,2$,  be
two totally regular sequences of spectral type of finite exponents
$s_{(i)}$, genus $p_{(i)}$, asymptotic domains
$AD_{\theta_{(i)},c_{(i)}}$, and orders $\alpha_{(i),N_{(i)}}\leq
0$. We assume in this section that the $\lambda_{(1),n_1}$ are real. However, note that all the results extend to the complex case with just minor modification in the definitions of the asymptotic domains. 
The sum sequence $S_{(0)}=S_{(1)}+S_{(2)}$ is a totally
regular sequence of spectral type with exponent
$s_{(0)}=s_{(1)}+s_{(2)}$, genus $p_{(0)}=[s_{(0)}]$,   asymptotic
domain $AD_{\theta_{(0)},c_{(0)}}$, where $\theta_{(0)}={\rm
max}(\theta_{(i)})$ and $c_{(0)}={\rm min}(c_{(i)})$, and order
$\alpha_{(0),N_{(0)}}={\rm min}(\alpha_{(i),N{(i)}})\leq 0$, by
the results at the end of Section \ref{s1}. The proof of the
following technical lemma is straightforward.

\begin{lem} \label{tec1} Assume $\{b_n\}_{n=1}^\infty $ is a sequence of
finite exponent and genus $p$, and $a$ is a real positive number, then
\[
\sum_{n=1}^\infty \sum_{j=1}^p \frac{(-1)^j}{j}a^j \left(
\frac{z^j}{(a+b_n)^j}-\frac{(1+z)^j}{b_n^j}+\frac{1}{b_n^j}\right)=
\sum_{k=1}^p E^p_k(a) z^k,
\]
where $E^p_k(a)=\sum_{n=1}^\infty\left(
\frac{(-1)^k}{k}\frac{a^k}{(a+b_n)^k}-\sum_{j=k}^p
\frac{(-1)^j}{j}\binom{j}{k}\frac{a^j}{b_n^j}\right)$.
\end{lem}

For example, if $p=k=1$,
\[
E_1^1(a)=\sum_{n=1}^\infty \left(-\frac{a}{a+b_n}+\frac{a}{b_n}\right)=\sum_{n=1}^\infty \frac{a^2}{b_n(a+b_n)}.
\]

\begin{teo} Suppose that $\alpha_{(1),N_{(1)}}<-p_{(2)}-1$,
and that $-\alpha_{(2),N_{(2)}}\geq s_{(1)}$. Then, the  sequence
$S_{(0)}=\{\lambda_{{(1)},n_1}+\lambda_{{(2)},n_2}\}_{n_i\in
\N_0}$ is spectrally decomposable over $S_{(1)}$, with power $1$, finite
length $\ell\leq N_{(2)}$, 
and asymptotic domain $AD_{\theta,c}$,
$\theta=\theta_{(0)}$ and $c={\rm
min}(c_{(0)},c_{(1)},1+c_{(2)}/\lambda_{(1),1})$. 
The length $\ell$ is the larger integer  
such that $-\alpha_{(2),\ell}\leq s_{(1)}$, where $\alpha_{(2),h}$ are the powers of the terms of the expansion of the $\Gamma$-function $\log\Gamma(-\lambda,S_{(2)})$.
\label{s6.t1}
\end{teo}
\begin{proof} As a double sequence, $S_{(0)}$ has finite exponents
$(s_{(0)}=s_{(1)}+s_{(2)}, s_{(1)},s_{(2)})$ and  genus
$(p_{(0)}=[s_{(0)}],p_{(1)},p_{(2)})$. We prove that $S_{(0)}$  is
spectrally decomposable with respect to $S_{(1)}$. We need to check
the requirements in Definition \ref{spdec}. We have that
$U=S_{(1)}$ is a totally regular sequence of spectral type of order
$u'=\alpha_{(1),N_{(1)}}$, with exponent $r_0=s_{(1)}$, genus
$q=p_{(1)}$ and domain $AD_{\theta_{(1)},c_{(1)}}$,
$0<c_{(1)}<\lambda_{(1),1}$. Let
$\widetilde{S}_{n_1}=\{\lambda_{(1),n_1}^{-1}\lambda_{n_1,n_2}\}_{n_2\in
\N_0}$, where $\lambda_{n_1,n_2}=\lambda_{{(1)},n_1}+\lambda_{{(2)},n_2}$. 
By Lemma \ref{lin1}, 
$\widetilde{S}_{n_1}$ is a
totally regular sequence of spectral type with exponent $s_{(2)}$,
genus $p_{(2)}$, spectral domain
$AD_{\theta_{(2)},1+c_{(2)}/\lambda_{(1),n_1}}$, and order
$\alpha_{(2),N_{(2)}}$, for all $n_1$. Thus we can take
$c'=1+c_{(2)}/\lambda_{(1),1}$ and $\theta'=\theta_{(0)}$ and the
first half of requirement (1) of Definition \ref{spdec} is
satisfied. Since $\alpha_{(i),N_{i)}}\leq 0$ by assumption, the second half of
the requirement (1) is also satisfied. In order to verify that the
requirement (2) is  satisfied we proceed as follows. We can decompose the $\Gamma$-function associated to $\tilde S_{n_1}$ as:
\begin{align*}
\log\Gamma(-\lambda,\widetilde{S}_{n_1})=-\sum_{n_2=1}^\infty&\left(\log
\left(1+\frac{(1-\lambda)
\lambda_{(1),n_1}}{\lambda_{(2),n_2}}\right) +\sum_{j=1}^{p_{(2)}}
\frac{(-1)^j}{j} \frac{((1-\lambda)
\lambda_{(1),n_1})^j}{\lambda_{(2),n_2}^j}\right.\\
&- \log \left(1+\frac{ \lambda_{(1),n_1}}{\lambda_{(2),n_2}}\right)\\
&-\sum_{j=1}^{p_{(2)}} \frac{(-1)^j}{j} \frac{(
\lambda_{(1),n_1})^j}{\lambda_{(2),n_2}^j} + \sum_{j=1}^{p_{(0)}}
\frac{(-1)^j}{j} \frac{(-\lambda
\lambda_{(1),n_1})^j}{(\lambda_{(1),n_1}+\lambda_{(2),n_2})^j}\\
&\left. -\sum_{j=1}^{p_{(2)}} \frac{(-1)^j}{j} \frac{((1-\lambda)
\lambda_{(1),n_1})^j}{\lambda_{(2),n_2}^j}+\sum_{j=1}^{p_{(2)}}
\frac{(-1)^j}{j} \frac{( \lambda_{(1),n_1})^j}{\lambda_{(2),n_2}^j}
\right)
\end{align*}
\begin{align*}
\log\Gamma(-\lambda,\widetilde{S}_{n_1})=&\log\Gamma((1-\lambda)\lambda_{(1),n_1},S_{(2)})-\log\Gamma(\lambda_{(1),n_1},S_{(2)})\\
&-\sum_{j=p_{(2)}+1}^{p_{(0)}} \frac{(-1)^j}{j} \sum_{n_2=1}^\infty
\frac{ \lambda_{(1),n_1}^j}{(\lambda_{(1),n_1}+\lambda_{(2),n_2})^j}
(-\lambda)^j \\
&-\sum_{j=1}^{p_{(2)}} \frac{(-1)^j}{j} \sum_{n_2=1}^\infty\hspace{-1pt}
\left(\frac{(-\lambda
\lambda_{(1),n_1})^j}{(\lambda_{(1),n_1}\hspace{-1pt}+\hspace{-1pt}\lambda_{(2),n_2})^j}
+\frac{(1-(1-\lambda)^j) \lambda_{(1),n_1}^j}{\lambda_{(2),n_2}^j} \right)\hspace{-2pt}.
\end{align*}

Applying Lemma \ref{tec1}, this gives
\beq\label{xx}
\begin{aligned}
\log\Gamma(-\lambda,\widetilde{S}_{n_1})
=&\log\Gamma((1-\lambda)\lambda_{(1),n_1},S_{(2)})-\log\Gamma(\lambda_{(1),n_1},S_{(2)})\\
&-\sum_{k=1}^{p_{(2)}}E^{p_{(2)}}_{k}(\lambda_{(1),n_1})
(-\lambda)^k,
\end{aligned}
\eeq
where
\[
E^{p_{(2)}}_{k}(\lambda_{(1),n_1})\hspace{-2pt}
=\hspace{-2pt}\left\{\begin{array}{ll}\hspace{-7pt}\sum_{n_2=1}^\infty\hspace{-2pt}\left(
\hspace{-2pt}\frac{(-1)^k}{k}\frac{\lambda_{(1),n_1}^k}{(\lambda_{(1),n_1}\hspace{-1pt}+\lambda_{(2),n_2})^k}\hspace{-2pt}
-\hspace{-2pt}\sum_{j=k}^{p_{(2)}}
\frac{(-1)^j}{j}\binom{j}{k}\frac{\lambda_{(1),n_1}^j}{\lambda_{(2),n_2}^j}\hspace{-1pt}\right)\hspace{-4pt},&\hspace{-4pt}
k\leq p_{(2)},\\
\hspace{-7pt}\frac{(-1)^k}{k} \sum_{n_2=1}^\infty \frac{
\lambda_{(1),n_1}^k}{(\lambda_{(1),n_1}+\lambda_{(2),n_2})^k},&\hspace{-4pt}
k> p_{(2)},\end{array}\right.
\]
are finite constants for each fixed $\lambda_{(1),n_1}$. Next,
recalling that $\frac{k}{j+k}\binom{j+k}{k}=\binom{j+k-1}{j}$, we can write
\[
E^{p_{(2)}}_k(-\lambda)=\frac{(-1)^{k+1}}{k}(-\lambda)^k\sum_{n_2=1}^\infty\hspace{0in}\left(
\frac{(-1)^{k+1}}{(\lambda-\lambda_{(2),n_2})^k}+\sum_{j=0}^{p_{(2)}-k}\hspace{-.1in}
\binom{j+k-1}{j}\frac{\lambda^j}{\lambda_{(2),n_2}^{j+k}}\right)\hspace{-2pt},
\]
and hence (compare with Lemmas  \ref{l0} and \ref{r2}) we find that
\beq\label{xx1}
\left.
E^{p_{(2)}}_k(\lambda_{(1),n_1})=\frac{(-1)^k}{k!}\left((-\lambda)^k
\frac{d^k}{d\lambda^k}
\log\Gamma(-\lambda,S_{(2)})\right)\right|_{\lambda=-\lambda_{(1),n_1}}.
\eeq

By hypothesis and Definition \ref{sss}, there exists an expansion of $\log\Gamma(-\lambda,S_{(2)})$ for large $\lambda$. This can be used to obtain the expansion of each of the terms appearing
in $\log\Gamma(-\lambda,\widetilde{S}_{n_1})$, as given in equation (\ref{xx}) (in particular also of the $E^{p_{(2)}}_k$, using equation (\ref{xx1})). After some computations and using Remark \ref{form}, we have the
expansion for large $\lambda_{(1),n_1}$ (where the $K_{k,h,i}$ are
some constants)
\begin{align}
\nonumber\log\Gamma(-\lambda,\widetilde{S}_{n_1})\hspace{-2pt} =&\sum_{h=0}^{N_{(2)}}
\left(a_{(2),\alpha_{(2),h},0}\left((1-\lambda)^{\alpha_{(2),h}}-1\right)
-\sum_{k=1}^{p_{(0)}} K_{k,h,0} (-\lambda)^k
\right)\lambda_{(1),n_1}^{\alpha_{(2),h}}\\
\label{ee1}
&+\sum_{l=0}^{p_{(2)}}\left(a_{(2),l,1}(1-\lambda)^{l}
\log(1-\lambda)-\sum_{k=1}^{p_{(0)}} K_{k,l,2} (-\lambda)^k
\right)      \lambda_{(1),n_1}^{l}\\
\nonumber&+\sum_{l=0}^{p_{(2)}}\left(a_{(2),l,1}\hspace{-2pt}\left((1-\lambda)^{l}-1\right)-\hspace{-2pt}\sum_{k=1}^{p_{(0)}}
K_{k,l,1} (-\lambda)^k
\right)\hspace{-2pt}\lambda_{(1),n_1}^{l}\hspace{-2pt}\log\lambda_{(1),n_1}\\
&\nonumber+o(\lambda_{(1),n_1}^{\alpha_{(2),N_{(2)}}}).
\end{align}

Comparing this expansion with the one required at point (2)
of Definition \ref{spdec}, observing that the polynomial appearing in the logarithmic terms vanish at $\lambda=0$, and observing that the following disequalities hold if the assumptions in the hypothesis are satisfied (note that in the present case ${\rm min}(-\alpha_{(2),0},-p_{(2)}-1)=-p_{(2)}-1$, see the beginning of Section \ref{s1}):
\begin{align*}
\alpha_{(1),N_{(1)}}&<-\alpha_{(2),0}< -\alpha_{(2),h}<
-\alpha_{(2),N_{(2)}},\\
\alpha_{(1),N_{(1)}}&<\alpha_{(1),0}\leq s_{(1)}\leq -\alpha_{(2),N_{(2)}},
\end{align*}
we deduce that all the requirements at point (2) of Definition \ref{spdec} are satisfied if we take as length of the decomposition the larger integer $\ell$ such that 
$-\alpha_{(2),\ell}\leq s_{(1)}$.
Note the value of $\kappa=1$. This concludes the proof of the theorem. In particular, note that 
equation (\ref{ee1}) gives the functions $\phi$,
as shown in the next corollaries. \end{proof}

Our next step is to find the explicit expression for the
functions $\phi_{\sigma_h}$ and the polynomials
$P_{\rho_l}$ that we need in order to apply Theorem \ref{sdl}. This is done comparing the formula in equation (\ref{ee1}) with the expansion of $\log\Gamma(-\lambda, \tilde S_n)$ given in point (2) of Definition \ref{spdec}.
For, it is convenient to analyze separately the three type of terms appearing in equation (\ref{ee1}). The terms of the first type are 
\[
\sum_{h=0}^{N_{(2)}}
\left(a_{(2),\alpha_{(2),h},0}\left((1-\lambda)^{\alpha_{(2),h}}-1\right)
-\sum_{k=1}^{p_{(0)}} K_{k,h,0} (-\lambda)^k
\right)\lambda_{(1),n_1}^{\alpha_{(2),h}}.
\]

They can be compared with the terms $\sum_{h=0}^N \phi_{\sigma_h}(\lambda)u_n^{-\sigma_h}$, 
and this gives the identification of the indices: $\sigma_h=-\alpha_{(2),h}$, for $h=0,\dots, H=N_{(2)}$.
The terms of the second type are
\[
\sum_{l=0}^{p_{(2)}}\left(a_{(2),l,1}(1-\lambda)^{l}
\log(1-\lambda)-\sum_{k=1}^{p_{(0)}} K_{k,l,2} (-\lambda)^k
\right)      \lambda_{(1),n_1}^{l}.
\]

Again, these are of the type of $\sum_{h=0}^N \phi_{\sigma_h}(\lambda)u_n^{-\sigma_h}$, 
and with $h=p_{(2)}-l$, we have: $\sigma_h=-l=-(p_{(2)}-h)$, for $h=0,\dots, H=p_{(2)}$. 
The terms of the last type are
\[
\sum_{l=0}^{p_{(2)}}\left(a_{(2),l,1}\hspace{-2pt}\left((1-\lambda)^{l}-1\right)-\hspace{-2pt}\sum_{k=1}^{p_{(0)}}
K_{k,l,1} (-\lambda)^k
\right)\hspace{-2pt}\lambda_{(1),n_1}^{l}\hspace{-2pt}\log\lambda_{(1),n_1},
\]
and give terms of the type $\sum_{l=0}^L P_{\rho_l}(\lambda)u_n^{-\sigma_h}\log u_n$. 
Here we identify $\rho_l=-(p_{(2)}-l)$, for $l=0,\dots, L=p_{(2)}$. Also note that the polynomial $P_{\rho_l}(\lambda)$ vanish at $\lambda=0$.
This analysis suggests to split the functions $\phi_{\sigma_h}(\lambda)$ appearing in point (2) of Definition \ref{spdec} as follows.

\begin{notation}\label{ass} We will write
\[
\sum_{h=0}^H \phi_{\sigma_h}(\lambda) u_n^{-\sigma_h}
=\sum_{h=0}^{\ell} \tilde\phi_{-\alpha_{(2),h}} (\lambda) u_n^{\alpha_{(2),h}}+
\sum_{l=0}^{p_{(2)}} \hat\phi_{-l}(\lambda) u_n^{l},
\]
with the identifications: $\sigma_h=-\alpha_{(2),h}$ in the first sum on the right side, and  $\sigma_h=h-p_{(2)}=-l$ and $\phi_{h-p_{(2)}}=\hat\phi_{-l}$  in the last sum. 
\end{notation}

Using this notation and the previous analysis, we prove the following corollary.

\begin{corol} \label{cor2}
\begin{align*}
 \tilde\phi_{-\alpha_{(2),h}}(\lambda)&=a_{(2),\alpha_{(2)},0}\left((1-\lambda)^{\alpha_{(2),h}}-1\right)
-\sum_{k=1}^{p_{(0)}} K_{k,h,0} (-\lambda)^k, &0\leq h\leq \ell,\\
 \hat\phi_{-l}(\lambda)&=a_{(2),l,1}(1-\lambda)^l\log(1-\lambda)
-\sum_{k=1}^{p_{(0)}} K_{k,l,2} (-\lambda)^k, &0\leq l\leq p_{(2)}.
\end{align*}
\end{corol}

\begin{corol} 
\begin{align*}
\tilde\Phi_{-\alpha_{(2),h}}(s)
&=a_{(2),\alpha_{(2),h},0}
\frac{\Gamma\left(s-\alpha_{(2),h}\right)}{s\Gamma\left(-\alpha_{(2),h}\right)} \hspace{.2in}(=0, ~{\rm if}~\alpha_{(2),h}\in \N),&0\leq h\leq \ell,\\
\hat\Phi_{-l}(s)&=(-1)^{l+1} l!a_{(2),l,1}
\frac{\Gamma(s-l)}{s},&0\leq l\leq p_{(2)}.
\end{align*}
\label{coroll}
\end{corol}
\begin{proof} We compute the function $\Phi_{\sigma_h}(s)$ using Corollary \ref{cor2} and 
the equations in the appendix. Note
that the second term in  both
$\tilde\phi_{\sigma_h}(\lambda)$ and $\hat\phi_{\sigma_h}(\lambda)$ gives no contribution, since the contour does not involve the
origin. Then, recalling the definition given in Lemma \ref{lp3}, we have
\[
\tilde\Phi_{-\alpha_{(2),h}}(s)=a_{(2),\alpha_{(2),h},0} \int_0^\infty
t^{s-1} \frac{1}{2\pi
i}\int_{\Lambda_{\theta,c}}\frac{\e^{-\lambda
t}}{-\lambda}\left((1-\lambda)^{\alpha_{(2),h}}-1\right)d\lambda
dt;
\]
the integral of the constant term vanishes, and we compute the first integral using equation (\ref{b4}) of the appendix. Note that the
expansion at $s=0$ is not uniform in the parameters. Therefore, we
must first compute the value of the function for each possible
allowed value of the parameters, and then obtain the expansion for
small $s$. By a similar argument using equation (\ref{b5}) of the appendix, we can also compute
\[
\hat\Phi_{-l}(s)=a_{(2),l,1} \int_0^\infty
t^{s-1} \frac{1}{2\pi
i}\int_{\Lambda_{\theta,c}}\frac{\e^{-\lambda
t}}{-\lambda}(1-\lambda)^l\log(1-\lambda)d\lambda
dt.
\]
\end{proof}

\begin{corol}\label{pippo} For $0\leq h\leq \ell$, $0\leq l\leq p_{(2)}$,  
\begin{align*}
\Rz_{s=0}\tilde\Phi_{-\alpha_{(2),h}}(s)&=
\left\{\begin{array}{ll}0,&\alpha_{(2),h}\in\N,\\
\psi\left(-\alpha_{(2),h}\right) a_{(2),\alpha_{(2),h},0},&\alpha_{(2),h}\not\in\N,
\end{array}\right.\\
\Ru_{s=0}\tilde\Phi_{-\alpha_{(2),h}}(s)&=
\left\{\begin{array}{ll}0,&\hspace{.4in}\alpha_{(2),h}\in\N,\\
a_{(2),\alpha_{(2),h},0},&\hspace{.4in}\alpha_{(2),h}\not\in\N,
\end{array}\right.\\
\Rd_{s=0}\tilde\Phi_{-\alpha_{(2),h}}(s)&=0,\\
\Rz_{s=0}\hat\Phi_{-l}(s)&=-k(l)a_{(2),l,1},\\
\Ru_{s=0}\hat\Phi_{-l}(s)&=-\psi(l+1)
a_{(2),l,1},\\
\Rd_{s=0}\hat\Phi_{-l}(s)&=-a_{(2),l,1},
\end{align*}
where $k(l)=\frac{(-1)^l}{2l!}\left(\frac{\pi^2}{3}+\psi^2(l+1)+\psi'(l+1)\right)$,
\end{corol}
\begin{proof} Using the expansion of $\Gamma(s+a)$ with $a\not\in -\N$ and $a\in-\N$ respectively, we
obtain, near $s=0$,
\[
\tilde\Phi_{-\alpha_{(2),h}}(s)= a_{(2),\alpha_{(2),h},0}\frac{1}{s}
+\psi\left(-\alpha_{(2),h}\right) a_{(2),\alpha_{(2),h},0} +O(s),
\]
and 
\[
\hat\Phi_{-l}(s)=-a_{(2),l,1}\frac{1}{s^2}-a_{(2),l,1}\frac{\psi(l+1)}{s} -a_{(2),l,1} k(l)+O(s).
\]
\end{proof}

The previous corollaries and the following lemma are the principal steps in the
proof of the main result of decomposition for the zeta function of
a sum of two sequences. In fact, they show that we can obtain all the
necessary information from the zeta functions of the summand sequences.

\begin{lem} 
\begin{align*}
A_{0,0}(s) &=\sum_{n_1=1}^\infty\left(
-\log\Gamma(\lambda_{(1),n_1},S_{(2)})+\sum_{j=0}^{p_{(2)}}a_{(2),j,1}\lambda_{(1),n_1}^j
\log\lambda_{(1),n_1}\right.\\
&\left.\hspace{.6in}+\sum_{h=0}^\ell
a_{(2),\alpha_{(2),h},0}\lambda_{(1),n_1}^{\alpha_{(2),h}}\right)\lambda_{(1),n_1}^{-s},\\
A_{0,1}(s)&=0.
\end{align*}
\label{laa}
\end{lem}
\begin{proof} We use Lemma \ref{lp2}. By definition, we need the asymptotic expansion of
$\log\Gamma(-\lambda,\widetilde{S}_{n_1})$ for large $\lambda$.
Such expansion can readily be obtained from equation (\ref{ee1}) at the
beginning of the proof of Theorem \ref{s6.t1}. We get
\begin{align*}
\log\Gamma(-\lambda,\widetilde{S}_{n_1})=&\log\Gamma((1-\lambda)\lambda_{(1),n_1},S_{(2)})-\log\Gamma(\lambda_{(1),n_1},S_{(2)})\\
&-\sum_{k=1}^{p_{(2)}}E^{p_{(2)}}_{k}(\lambda_{(1),n_1})
(-\lambda)^k-\sum_{k=p_{(2)}+1}^{p_{(0)}}
E^{p_{(2)}}_{k}(\lambda_{(1),n_1}) (-\lambda)^k,
\end{align*}
where the $E^{p_{(2)}}_{k}(\lambda_{(1),n_1})$ are finite constants.
We just need the coefficients of the constant and of the logarithmic
terms, so there are not contributions from the last two terms. Using the expansion given in Remark \ref{form} for the first term, we obtain
\begin{align*}
\log\Gamma(-\lambda,\widetilde{S}_{n_1})=&-\log\Gamma(\lambda_{(1),n_1},S_{(2)})+\sum_{j=0}^{p_{(2)}}a_{(2),j,1}\lambda_{(1),n_1}^j (1-\lambda)^j
\log\lambda_{(1),n_1}\\
&+\sum_{j=0}^{p_{(2)}}a_{(2),j,1}\lambda_{(1),n_1}^j (1-\lambda)^j
\log(-\lambda)\\
&+\sum_{h=0}^{N_{(2)}} a_{(2),\alpha_{(2),h},0}
((1-\lambda)\lambda_{(1),n_1})^{\alpha_{(2),h}}\\ 
&+\sum_{j=0}^{p_{(2)}}a_{(2),j,1}\lambda_{(1),n_1}^j
(1-\lambda)^j \log\left(1+\frac{1}{-\lambda}\right)
+o((-\lambda)^{\alpha_{(2),N_{(2)}}}).
\end{align*}

Thus the coefficient of the logarithmic term is $a_{0,1,n_1}=\sum_{j=0}^{p_{(2)}}a_{(2),j,1}\lambda_{(1),n_1}^j$. 
To get the one of the constant term is a bit harder. From one side, expanding
\[
(1-\lambda)^j \log\left(1+\frac{1}{-\lambda}\right)=\sum_{m=0}^j
\binom{j}{m} (-\lambda)^m \sum_{k=1}^\infty
\frac{(-1)^{k+1}}{k}(-\lambda)^{-k},
\]
we see that the constant term of this part is 
$\sum_{m=1}^j \binom{j}{m}\frac{(-1)^{m+1}}{m}=\sum_{m=1}^j
\frac{1}{m}$.

On the other side, recalling the powers of the binomial, we see
that the unique contributions to the constant term of the term
$\sum_{h=0}^{N_{(2)}} a_{(2),\alpha_{(2),h},0}
((1-\lambda)\lambda_{(1),n_1})^{\alpha_{(2),h}}$
in the expansion above are by the terms with positive integer
exponents, namely
\[
\sum_{h=0,\alpha_{(2),h}\in \N}^{N_{(2)}} a_{(2),\alpha_{(2),h},0}
\lambda_{(1),n_1}^{\alpha_{(2),h}}= 
\sum_{h=0,\alpha_{(2),h}\in \N}^{\ell} a_{(2),\alpha_{(2),h},0}
\lambda_{(1),n_1}^{\alpha_{(2),h}}.
\]

The reason for the equality is the following. First, notice that the sum over $h$ is a finite sum even when $N_{(2)}=\infty$. In fact, since the $\alpha_{(2),h}$ form a decreasing sequence with unique accumulation point at $-\infty$, there can be at most a finite numbers of non negative entries. Second, since $s_{(1)}$ is non negative, by the same reasoning as in the proof of Theorem \ref{s6.t1}, we can restrict the sum to $h\leq \ell$.

Collecting all together, the complete constant term is
\begin{align*}
a_{0,0,n_1}=&-\log\Gamma(-\lambda_{(1),n_1},S_{(2)})+\sum_{j=0}^{p_{(2)}}a_{(2),j,1}\lambda_{(1),n_1}^j
\log\lambda_{(1),n_1}\\
&+\sum_{j=1}^{p_{(2)}}a_{(2),j,1}\sum_{m=1}^j
\frac{1}{m}\lambda_{(1),n_1}^j+ \sum_{h=0, \alpha_{(2),h}\in\N}^\ell
a_{(2),\alpha_{(2),h},0} \lambda_{(1),n_1}^{\alpha_{(2),h}}.
\end{align*}

The last contribution to the $A_{x,j}$ as given in Lemma \ref{lp2} comes from the coefficients in the expansion for
large $\lambda$ of the functions $\phi_{\sigma_h}(\lambda)$. 
Actually, we only need the coefficients $b_{\sigma_h,0,k}$. Using 
the formulas given in Corollary \ref{cor2}, and comparing with the formulas in Lemma \ref{lll1}, we obtain after some calculations
\[
\begin{array}{lll}\tilde b_{-\alpha_{(2),h},0,1}=0, 
&
\tilde b_{-\alpha_{(2),h},0,0}=\left\{\begin{array}{ll}-a_{(2),\alpha_{(2),h},0},&\alpha_{(2),h}\not\in\N,\\
0&\alpha_{(2),h}\in\N,
\end{array}\right.&0\leq h\leq \ell,\\
\hat b_{-l,0,1}=a_{(2),l,1},&
\hat b_{-l,0,0}=\left\{\begin{array}{ll}\sum_{m=1}^l \frac{1}{m} a_{(2),l,1}&l\in\N_0,
\\0&l=0.
\end{array}\right.&0\leq l\leq p_{(2)}.
\end{array}
\]

Rearranging, we have
\begin{align*}
b_{\sigma_h,0,1}&=\left\{
\begin{array}{ll}0&\sigma_h=-l\not\in-\N,\\a_{(2),l,1}& \sigma_h=-l\in-\N,\end{array}\right.\\
b_{\sigma_h,0,0}&=\left\{
\begin{array}{ll}-a_{(2),\alpha_{(2),h},0}&0\leq h\leq N,~\sigma_h=-\alpha_{(2),h}\not\in-\N,\\
\sum_{m=1}^l \frac{1}{m} a_{(2),l,1}&
0\leq l\leq p_{(2)},~\sigma_h=-l\in-\N_0,\\0&\sigma_h=-l=0.\end{array}\right.
\end{align*}

Substituting these values in the formulas of Lemma \ref{lp2}
we have the thesis.
\end{proof}

We give now the first terms of the expansion at $s=0$ of $\zeta(s,S_{(0)}=S_{(1)}+S_{(2)})$. Here the coefficients $a_{(i),x,h}$ and $c_{(i),x,h}$ are as in Remark \ref{form}.

\begin{teo} \label{tt1} Suppose that $\alpha_{(1),N_{(1)}}<-p_{(2)}-1$, and $
-\alpha_{(2),N_{(2)}}\geq s_{(1)}$. Then, the zeta function associated to the sum sequence 
$S_{(0)}=\{\lambda_{n_1,n_2}=\lambda_{{(1)},n_1}+\lambda_{{(2)},n_2}\}_{n_i\in
\N_0}$ is regular at $s=0$ and 
\begin{align*}
&\zeta(0,S_{(0)})=a_{(1),0,1}a_{(2),0,1}\hspace{-1.5pt}+\hspace{-1.5pt}\sum_{j=1}^{p_{(2)}} (-1)^{j+1}j a_{(1),-j,0}a_{(2),j,1}
\hspace{-1.5pt}+\hspace{-1.5pt}\sum_{h=0}^\ell \hspace{-1pt}\frac{
a_{(1),-\alpha_{(2),h},0}a_{(2),\alpha_{(2),h},0}}
{\Gamma(-\alpha_{(2),h})\Gamma(\alpha_{(2),h})}\\
&\hspace{40pt}= \sum_{h=0}^{\ell}
c_{(1),-\alpha_{(2),h},0}c_{(2),\alpha_{(2),h},0},\\
&\zeta'(0,S_{(0)})=-\sum_{l=0}^{p_{(2)}}a_{(2),l,1}\left(
\zeta'(-l,S_{(1)})+(\gamma+\psi(l+1))\zeta(-l,S_{(1)})\right)\\
&+\hspace{-8pt}\sum_{h=0,\alpha_{(2),h}\not\in\N}^\ell
\hspace{-8pt}a_{(2),\alpha_{(2),h},0}
\left(\Rz_{s=-\alpha_{(2),h}}\zeta(s,S_{(1)})+
(\gamma\hspace{-1pt}+\hspace{-.7pt}\psi(-\alpha_{(2),h}))\Ru_{s=-\alpha_{(2),h}}\zeta(s,S_{(1)})\right)\\
&+\log \prod_{n_1=1}^\infty\e^{-\sum_{j=0}^{p_{(2)}}a_{(2),j,1}\lambda_{(1),n_1}^j
\log\lambda_{(1),n_1}-\sum_{h=0}^\ell a_{(2),\alpha_{(2),h},0}\lambda_{(1),n_1}^{\alpha_{(2),h}}}
\Gamma(\lambda_{(1),n_1},S_{(2)}).
\end{align*}
\end{teo}
\begin{proof} By Theorem \ref{s6.t1}, $S_{(0)}$ is spectrally decomposable
over $S_{(1)}$ with length $\ell$ and power $\kappa=1$. Hence, we apply Theorem \ref{sdl}. By  Corollary \ref{pippo}, 
$\Res{2}{0}\tilde\Phi_{-\alpha_{(2),h}}(s)=0$, and by Proposition \ref{aaab3}, $\Res{1}{-l}\zeta(s,S_{(1)})=0$, 
therefore $\zeta(s,S_{(0)})$ is regular at $s=0$. Next, we use the formulas of Theorem \ref{sdl} to compute
$\zeta(0,S_{(0)})$ and $\zeta'(0,S_{(0)})$. Using the notation introduced in \ref{ass}, Corollary \ref{pippo}, Lemma \ref{laa} and Proposition \ref{aaab3}, we have
\begin{align*}
\Rz_{s=0}\zeta(s,S)
=&-\sum_{j=0}^{p_{(2)}} a_{(2),j,1} \zeta(-j,S_{(1)}) +\sum_{h=0}^\ell
\Ru_{s=-\alpha_{(2),h}}\zeta(s,S_{(1)})a_{(2),\alpha_{(2),h},0},
\end{align*}

A further application of Proposition \ref{aaab3} gives the first formula stated in the theorem for $\zeta(0,S)$. The second formula in the theorem follows applying Proposition \ref{z3} (see also the remarks below). The proceedure followed in order to obtain the formula for the derivative is similar but much longer, and we omit the details.
\end{proof}

It is important to observe that, applying Propositions \ref{z3} and \ref{aaab3}, the formulas given in Theorem \ref{tt1}, can be written using only particular values or residues, i.e. zeta invariants, of the zeta functions associated to the two sequences $S_{(i)}$. Explicit formulas are given in the next corollary. 

\begin{corol}\label{formulas}
\begin{align*}
\zeta(0,S_{(0)})=&\zeta(0,S_{(1)})\zeta(0,S_{(2)})+ \sum_{j=1}^{p_{(2)}} \frac{(-1)^j}{j}
\zeta(-j,S_{(1)})\Ru_{s=j}\zeta(s,S_{(2)})\\
&+\sum_{l=1}^{p_{(1)}} \frac{(-1)^l}{l} \zeta(-l,S_{(2)})\Ru_{s=l}\zeta(s,S_{(1)})\\
&+\hspace{-2pt}\sum_{h=0,\alpha_{(2),h}\not\in\Z}^{\ell}\hspace{-2pt}\Gamma(\alpha_{(2),h})
\Gamma(-\alpha_{(2),h})\Ru_{s=-\alpha_{(2),h}}\zeta(s,S_{(1)})
\Ru_{s=\alpha_{(2),h}}\zeta(s,S_{(2)}),
\end{align*}
\begin{align*}
\zeta'(0,S_{(0)})
=&\zeta(0,S_{(2)})\zeta'(0,S_{(1)})\hspace{-.6pt}+\hspace{-2pt}\sum_{h=0,\alpha_{(2),h}\not\in\Z}^\ell\hspace{-2pt} \Gamma(\alpha_{(2),h})\Gamma(-\alpha_{(2),h})
\Ru_{s=\alpha_{(2),h}}\hspace{-2pt}\zeta(s,S_{(2)})
\\
&\hspace{38pt}\left(\Rz_{s=-\alpha_{(2),h}}\zeta(s,S_{(1)})+
(\gamma+\psi(-\alpha_{(2),h}))\Ru_{s=-\alpha_{(2),h}}\zeta(s,S_{(1)})\right)\\
&+\sum_{l=1}^{p_{(1)}} \frac{(-1)^{l}}{l}\zeta(-l,S_{(2)})\hspace{-3pt}\left(\hspace{-1.5pt}\Rz_{s=l}\zeta(s,S_{(1)})+
(\gamma+\psi(l))\Ru_{s=l}\zeta(s,S_{(1)})\hspace{-1.5pt}\right)\\
&+\sum_{l=1}^{p_{(2)}}\frac{(-1)^{l}}{l}\Ru_{s=l}\zeta(s,S_{(2)})\hspace{-.6pt}\left(
\zeta'(-l,S_{(1)})+(\gamma+\psi(l\hspace{-2pt}+\hspace{-2pt}1))\zeta(-l,S_{(1)})\right)\\
&-\log \prod_{n_1=1}^\infty\left(
\e^{-\sum_{l=1}^{p_{(2)}}\frac{(-1)^{l}}{l}\Rs{1}{l}\zeta(s,S_{(2)})\lambda_{(1),n_1}^l
\log\lambda_{(1),n_1}}
\right.\\
&\e^{\sum_{h=0,\alpha_{(2),h}\not\in\Z}^\ell
\Gamma(\alpha_{(2),h})\Gamma(-\alpha_{(2),h})
\Rs{1}{\alpha_{(2),h}}\hspace{-.2in}\zeta(s,S_{(2)})\lambda_{(1),n_1}^{\alpha_{(2),h}}}\\
&\e^{\sum_{l=1}^{p_{(1)}} \frac{(-1)^{l}}{l}\zeta(-l,S_{(2)})\lambda_{(1),n_1}^{-l}}
\e^{-\sum_{l=1}^{p_{(2)}} \frac{(-1)^{l}}{l}\left(\Rs{0}{l}\zeta(s,S_{(2)})-\frac{1}{l}\Rs{1}{l}\zeta(s,S_{(2)})\right)\lambda_{(1),n_1}^l}\\
&\left.
\e^{-\zeta(0,S_{(2)})\log\lambda_{(1),n_1}-\zeta'(0,S_{(2)})}
\prod_{n_2=1}^\infty
\left(1+\frac{\lambda_{(1),n_1}}{\lambda_{(2),n_2}}\right)
\e^{\sum_{j=1}^{p_{(2)}}\frac{(-1)^j}{j}\frac{\lambda_{(1),n_1}^j}{\lambda_{(2),n_2}^j}}\right)\hspace{-2pt}.
\end{align*}
\end{corol}

The formulas in Theorem \ref{tt1} and its corollary look quite complicate and present an overweight in notation. This is a natural problem in very general zeta regularization formalisms. However, these formulas contain in fact more information with respect to the ones presented in Theorem \ref{sdl}, and this justifies the work that was necessary to obtain them. We clarify this point in the following remarks, and with two applications in the following Section \ref{s4}.
First, note the simple formula in the theorem  for the value of $\zeta(0,S_{(0)})$, where only the coefficients $c_{(i),x,h}$ of the asymptotic expansions of the heat functions of the two sequences (see Proposition \ref{coeff}) appear. This is a consequence of the fact that the value of $\zeta(s,S_{(0)})$ at zero is local, and hence depends only on the expansion of the heat function $f(t,S_{(0)})$ (by Mellin transform). Second, as in Theorem \ref{spdec},  in all the formulas there are two types of contributions (see Lemmas \ref{lp1} and \ref{lp2}): those coming from the regular part, $\zeta_1$ of Lemma  \ref{lp3},  and those coming from the singular part, $\zeta_2$ and $\zeta_3$. However, in the present case is characterized by two remarkable properties: first, $s=0$ is a regular point, and second, there are no regular contribution to $\zeta(0,S_{(0)})$. Both properties are consequences of the fact that $S_{(0)}$ is defined as linear combination of two sequences of spectral type, and hence the coefficient $a_{0,1,n_1}$ of the logarithmic term of $\log\Gamma(-\lambda, \tilde S_{n_1})$  coincides with the sum of the coefficients $b_{\sigma_h,0,1}$ of the logarithmic terms of the functions $\phi_{\sigma_h}(\lambda)$ (see the proof of Lemma \ref{laa}). This gives $A_{0,1}(s)=0$. 
Third, by similar analysis, we distinguish in $\zeta'(0,S_{(0)})$ two contributions. The singular one, coming from $\zeta_2$ and $\zeta_3$, is given by the first two lines in the formula in the theorem, and consists of two finite sums of products of the coefficients $a_{(i),x,j}$ of the expansions of the logarithmic Gamma function of the two sequences  (see Lemma \ref{lll1}) and/or particular values of the two zeta functions $\zeta(s,S_{(i)})$ (the conversion is done using Proposition \ref{aaab3}). This contribution decomposes completely on information coming from the zeta functions of the two summand sequences. The regular contribution has a more complicate description. It is the one in the last line of the theorem, described more explicitely in the last four lines of the corollary, and has the following analytic interpretation: it is a 'regularized' double Weierstrass canonical product, in the following sense. The inner product $\prod_{n_2=1}^\infty$ defines a proper Weierstrass product. 
Next, consider the sequence of functions $\log\Gamma(\lambda_{(1),n_1},S_{(2)})$. 
Then, the sum $\sum_{n_1=1}^\infty \lambda_{(1),n_1}^{-s}\log\Gamma(\lambda_{(1),n_1},S_{(2)})$ is well defined provided that we subtract all the divergent terms, namely all the terms in the expansion for large $\lambda_{(1),n_1}$ of the logarithmic Gamma function that behave like $n_1^x$ with $x\geq-1$. This is where the outer  product comes in, and the subtracted terms are those in the exponents. They correspond  bijectively to the terms appearing in the singular part. 



\section{Two applications}
\label{s4}

\subsection{The Kronecker first limit formula for sequences of spectral type}
\label{s4.a}

In this section we give an application to the study of  some special functions appearing in number theory. More precisely, we show how a classical formula of analytic number theory, the Kronecker first limit formula, extends when the sequence of the positive integer numbers is replaced by any totally regular sequence of spectral type, and the Dedeckind eta function appearing in that formula is replaced by a generalization that satisfies the same functional equation. As observed in the introduction, this example is presented in a quite detailed way, with the aim of helping the reader in understanding how to apply the results of the previous sections. 

Recall that the Dedekind  zeta function of a number field $\K$ is defined by
\[
\zeta_\K(s)=\sum_{\mathsf{a}} n(\mathsf{a})^{-s},
\]
where $\mathsf{a}$ varies over the integral ideals of $\K$ and $n(\mathsf{a})$ denotes their absolute norm. The Dedekind zeta function has an analytic continuation to the whole complex plane with a simple pole at $s=1$, satisfies a functional equation, and  decomposes into a finite sum $\zeta_\K(s)=\sum_{\A}\zeta_\A(s)$, where $\A$ runs over the ideal class group of $\K$. The Laurent expansion at $s=1$ of $\zeta_\A(s)$, is called the {\it Kronecker first limit formula in the number field} $\K$, and the residue is independent of the ideal class $\A$ chosen. This fact is at the basis of the analytic determination of the class number of $\K$. The evaluation of the finite part $\rho(\A)$ is great importance in analytic number theory. Using the functional equation of $\zeta_\A(s)$, we can reduce the Kronecker first limit formula to the expansion at $s=0$
\beq\label{e1}
\zeta_\A(s)=\zeta_\A(0)+\zeta'_\A(0) s+O(s^2),
\eeq
that shows the importance of the evaluation of the derivative of the zeta function in this context.
The {\it classical Kronecker first limit formula} is the particular instance of equation (\ref{e1}) when the field is an imaginary quadratic field. In  this case, the zeta function $\zeta_\A(s)$ reduces to the following zeta function (up to some known factors), called {\it Eisenstein (or Epstein) zeta function} 
\[
E(\tau,s)=y^{s}{\sum_{m,n\in\Z}}^{\hspace{-6pt}\prime} |m\tau+n|^{-2s},
\]
where $\tau=x+iy$ with $y>0$. Consider for simplicity the case $x=0$. Then, $E(\tau,s)=y^s\zeta(s,y)$, where
\beq\label{zt}
\begin{aligned}
\zeta(s,0,y)=&{\sum_{m,n\in\Z}}^{\hspace{-6pt}\prime} (y^2m^2+n^2)^{-s}\\
=&2y^{-2s}\zeta_R(2s)+2\zeta_R(2s)+4\zeta(s,S+y^2S),
\end{aligned}
\eeq
where $S$ is the sequence $\{n^2\}_{n=1}^\infty$. This is a totally regular sequence of spectral type of exponent $\ec=\frac{1}{2}$, genus $\ge=0$, infinite order $N=\infty$, $\alpha_N=-\infty$. For we have asymptotic expansions both for the associated logarithmic Gamma function and for the associated heat function (see \cite{Spr4} Section 3.1):
\begin{align*}
\log\Gamma(-\lambda,S)=&-\log\prod_{n=1}^\infty\left(1-\frac{\lambda}{n^2}\right)=-\log\frac{{\rm sh}\pi\sqrt{-\lambda}}{\pi\sqrt{-\lambda}}\\
=&-\pi\sqrt{-\lambda}+\frac{1}{2}\log(-\lambda)+\log 2\pi+O\left(\e^{-2\pi\sqrt{-\lambda}}\right),\\
f(t,S)-1=&\sum_{n=1}^\infty \e^{-n^2 t}=\frac{\sqrt{\pi}}{2}t^{-\frac{1}{2}}-\frac{1}{2}+O\left(\e^{-\frac{1}{t}}\right).
\end{align*}

From Remark \ref{form} we obtain the indices $\alpha_j=\frac{1-j}{2}$, $j=0,\infty$, and
\[
\left\{\begin{array}{l}\alpha_0=\frac{1}{2}\\ \alpha_1=0\\ \alpha_2=-\frac{1}{2},\end{array}\right.
\]
and the coefficients
\[
\left\{\begin{array}{ll}c_{\alpha_0,0}=c_{\frac{1}{2},0}=\frac{\sqrt{\pi}}{2}&\\
c_{\alpha_1,0}=c_{0,0}=-\frac{1}{2}&\\
c_{\alpha_2,0}=c_{-\frac{1}{2},0}=0&\\
c_{\alpha_j,0}=c_{\frac{1-j}{2},0}=0&j>1.\\
\end{array}\right.
\]

Thus $S$ is totally regular. Moreover (compare with Proposition \ref{z3})
\[
\left\{\begin{array}{ll}a_{\alpha_0,0}=a_{\frac{1}{2},0}=\Gamma(-\alpha_0)c_{\alpha_0,0}
=\Gamma(-\frac{1}{2})c_{\frac{1}{2},0}=-\pi&\\
a_{\alpha_1,0}=a_{0,0}=-\zeta'(0,S)&\\
a_{\alpha_j,0}=a_{\frac{1-j}{2},0}=0&j>1,\\
\end{array}\right.
\]
\[
\left\{\begin{array}{ll}a_{\alpha_0,1}=a_{\frac{1}{2},1}=0&\\
a_{\alpha_1,1}=a_{0,1}=-\zeta(0,S)=-c_{0,0}=\frac{1}{2}&\\
a_{\alpha_j,1}=a_{\frac{1-j}{2},1}=0&j>1,\\
\end{array}\right.
\]

The associate zeta function is $\zeta(s,S)=\zeta_R(2s)$, and hence: $\zeta(0,S)=\zeta_R(0)=-\frac{1}{2}$, $\zeta'(0,S)=2\zeta_R'(0)=-\log 2\pi$.

Next, consider the sum sequence $S_{(0)}=y^2S+S=S_{(1)}+S_{(2)}$. We have  $s_{(i)}=\frac{1}{2}$, $p_{(i)}=0$, so that $s_{(0)}=1$ and $p_{(0)}=1$. Since $-\infty=\alpha_{(1),N_{(1)}}<-p_{(2)}-1=-1$, and $\infty=-\alpha_{(2),N_{(2)}}> s_{(1)}=\frac{1}{2}$, it follows by Theorem \ref{s6.t1} that $S_{(0)}$ is spectrally decomposable over $S_{(1)}$. We deduce the length of the decomposition. Since:
\[
\begin{matrix}h&=&0&1&2&3&4&\dots\\
-\alpha_{(2),h}&=&-\frac{1}{2}&0&\frac{1}{2}&1&\frac{3}{2}&\dots
\end{matrix}
\]
so, in order to satisfy the condition: $\ell$ is the large integer such that $-\alpha_{(2),\ell}\leq s_{(1)}=\frac{1}{2}$,
we obtain $\ell=2$, and $[|\alpha_{(2),\ell}|]=[|\alpha_{(2),2}|]=0$.
Proceeding as above, we easily obtain the following information on the sequence $S_{(1)}=y^2S$. 
\[
f(t,y^2S)-1=\sum_{n=1}^\infty \e^{-y^2n^2 t}=\frac{\sqrt{\pi}}{2y}t^{-\frac{1}{2}}-\frac{1}{2}+O(\e^{-\frac{1}{t}}),
\]
so the unique coefficient that changes is $c_{\frac{1}{2},0}=\frac{\sqrt{\pi}}{2y}$. Also, $\zeta(s,y^2S)=y^{-2s}\zeta_R(2s)$, gives $\zeta(0,y^2S)=\zeta_R(0)=-\frac{1}{2}$, $\zeta'(0,y^2S)=-2\zeta_R(0)\log y+2\zeta_R'(0)=\log y-\log 2\pi$, $\Rz_{s=-\frac{1}{2}}\zeta(s,y^2S)=\zeta(-\frac{1}{2},y^2S)=y\zeta_R(-1)=-\frac{y}{12}$.
We can now apply Corollary \ref{formulas}. After some calculations, we obtain
\begin{align*}
\zeta(0,y^2 S+S)&=\frac{1}{4},\\
\zeta'(0,y^2 S+S)
&=-\zeta(0,y^2S+S)\log y^2-\log\eta(iy,S),
\end{align*}
where
\[
\eta(iy,S)=\frac{1}{\sqrt{2\pi}}\e^{\frac{\pi}{12}y}\prod_{n=1}^\infty \left(1-\e^{-2\pi yn}\right),
\]
and the {\it Kronecker first limit formula} for the zeta function $\zeta(s,y)$ is:
\begin{align}\label{k1}
\zeta(s,y)
&=-1-2\left(\log2\pi+2\log|\eta_D(iy)|\right)s+O(s^2).
\end{align}

The function $\eta(\tau,S)$ appearing  the limit formula is (up to a constant) the well known  {\it Dedekind eta function}: $\eta(\tau,S)=\frac{1}{2\pi}\eta_D(\tau)$, that satisfied to the functional equation: 
\beq\label{feta}
\eta_D\left(-\frac{1}{\tau}\right))=\sqrt\frac{\tau}{i}\eta_D(\tau).
\eeq

\begin{rem}\label{rtorus} Note that the zeta function $\zeta(s,yS+S)$  gives the zeta function associated to the Laplace operator $\Delta_{T_y}$ on a torus $S^1\times S^1_y$ (where $S^1_y$ is the circle of radius $y$) (compare with \cite{OS} Section 5).
\end{rem}

Consider now the case where the sum sequence is $S_{(0)}=y S+S$, with $S$ some generic totally regular sequence of spectral type and $y$ some positive real number. Then, $\ell=N_{(1)}=N_{(2)}=N$, $\alpha_{(1),h}=\alpha_{(2),h}=\alpha_h$, $p_{(1)}=p_{(2)}=\ge$, and by Theorem \ref{s6.t1}, $yS+S$ is spectrally decomposable over $yS$ with finite length if $\alpha_N<-\ge-1$. Assuming that this is the case, we can perform the same analysis as above, and we obtain all necessary information on the indices, the coefficients, the length etc. A tedious but straightforward application of the formulas of Corollary \ref{formulas} gives 
\[
\zeta'(0,S+yS)=-\zeta(0,S+yS)\log y-\log\eta(i\sqrt{y},S),
\]
where we define the {\it generalized Dedekind eta function} by
\begin{align*}
\eta(i\sqrt{y},S)=&
\e^{\sum_{l=0}^{\ge}a_{l,1}\left(
\zeta'(-l,S)+(\gamma+\psi(l+1))\zeta(-l,S)\right)y^l}\\
&\e^{-\sum_{h=0,\alpha_{h}\not\in\N}^\ell
a_{\alpha_{h},0}\left(\Rs{0}{-\alpha_h}\hspace{-10pt}\zeta(s,S)+
(\gamma\hspace{-1pt}+\hspace{-.7pt}\psi(-\alpha_{h}))\Rs{1}{-\alpha_h}\hspace{-10pt}
\zeta(s,S)\right)y^{\alpha_h}}\\
\prod_{n_1=1}^\infty&\e^{\sum_{j=0}^{\ge}a_{j,1}\lambda_{n_1}^j y^j
\log y \lambda_{n_1}+\sum_{h=0}^\ell a_{\alpha_{h},0}\lambda_{n_1}^{\alpha_{h}}y^{\alpha_{h}}}
\prod_{n_2}^\infty \left(1+\frac{\lambda_{n_1}}{\lambda_{n_2}}y\right)\e^{\sum_{j=1}^\ge \frac{(-y)^j}{j}\frac{\lambda_{n_1}^j}{\lambda_{n_2}^j}}\hspace{-2pt},
\end{align*}
and reduces to the classical Dedekind eta function when $S=\{n^2\}_{n=1}^\infty$.

Summing up, we get the {\it generalized Kronecker first limit formula} for the zeta function $\zeta(s,S+yS)$: 
\beq
\label{eta0}
\zeta(s,S+yS)=\zeta(0,S+yS)-\left(\zeta(0,S+yS)\log y+\log\eta(i\sqrt{y},S)\right)s +O(s^2).
\eeq

The formula in equation (\ref{eta0}) is a striking result since, for any totally regular $S$-sequence $S$, the new spectral function $\eta(\tau,S)$ satisfies precisely the same functional equation (\ref{feta}) satisfied by the classical Dedekind eta function appearing in the classical Kronecker first limit formula (\ref{k1}). In fact, the symmetry of the double zeta function under twist of summation indices gives
\[
\zeta(s,S+yS)=y^{-s}\zeta(s,S+\frac{1}{y}S), 
\]
and this implies that
\begin{align*}
\zeta(0,S+yS)=&\zeta(0,S+\frac{1}{y}S),\\
\zeta'(0,S+yS)=&-\zeta(0,S+\frac{1}{y}S)\log y +\zeta'(0,S+\frac{1}{y}S).
\end{align*}

Substitution in equation (\ref{eta0}) gives the following  {\it functional equation} for the generalized eta function:
\[
\eta\left(\frac{i}{y},S\right)=y^{2\zeta(0,S+y^2 S)}\eta(iy,S),
\]
to be compared with equation (\ref{feta}).

\subsection{Zeta determinant of a product space}\label{s4.2} For a compact connected Riemannian manifold $M$,
if $Sp_+\Delta_M=\{\lambda_n\}$ is the positive spectrum of the metric Laplacian $\Delta_M$, the associate zeta function is defined by \cite{RS}
\[
\zeta(s,Sp_+\Delta_M)=\sum \lambda_n^{-s},
\]
and the zeta determinant by
\[
\det_\zeta \Delta_M=\e^{-\zeta'(0,Sp_+ \Delta_M)}.
\]

In this section we answer the following question: given
two compact connected Riemannian manifolds $M_{(1)}$ and
$M_{(2)}$, what is the relation between the zeta
determinant of the Laplace operator on the product space ${\rm
det}_\zeta(\Delta_{M_{(1)}\times M_{(2)}})$ and those of the
factors ${\rm det}_\zeta(\Delta_{M_{(i)}})$? In order to answer this question, we first recall some basic facts about spectral properties  of operators on manifolds, zeta determinants and analytic properties of the other spectral functions. In particular, we state these well known properties in the language of spectral sequences of Section \ref{s1}.

Let $M$ be a compact connected Riemannian manifold of dimension
$m$. The (negative of the) metric Laplacian $\Delta_M$ has a positive real spectrum $0<Sp_+\Delta_M=\{\lambda_1\leq\lambda_{2}\leq \dots\}$ (plus 
the possible null eigenvalue  if $M$ has no boundary), whose behavior for large $n$ is known by the Weyl formula (see for example \cite{Gil}).
There is a full asymptotic expansion for the trace of the heat operator for small $t$, 
\[
{\rm Tr}_{L^2}\e^{-t\Delta_M}=t^{-\frac{m}{2}}\sum_{j=0}^\infty e_j t^\frac{j}{2},
\]
where the coefficients  depend only on local
invariants constructed from the metric tensor, and are in
principle calculable from it (and all the coefficients of odd index vanish
if the manifold has no boundary).
We can apply the method of Section \ref{s1} to the sequence $Sp_+ \Delta_M$.
Using Propositions \ref{coeff}, \ref{zeta1}, \ref{z3}, \ref{aaab3} and Theorem \ref{zeta2}, we can prove the following results.

\begin{prop} \label{s3.l1} The sequence $Sp_+ \Delta_M$
of the eigenvalues of the metric Laplacian
on a compact
connected Riemannian manifold of dimension $m$, is a totally regular sequence of spectral type, with finite
exponent $\ec=\frac{m}{2}$, genus 
\[
\ge=[\ec]=\left\{\begin{array}{cc}\frac{m}{2}&m ~{\rm even},\\
\frac{m-1}{2}&m ~{\rm odd},\end{array}\right.
\]
spectral sector $\Sigma_{\epsilon,c}$, with  $0< c< \lambda_1$, asymptotic
domain $AD_{\epsilon,c}$, and infinite order.
\end{prop}

We have the following formulas for the coefficients in the expansion of the heat function
$\alpha_h=\frac{m-h}{2}$, $c_{\frac{m-h}{2},1}=0$,  
\[
c_{\alpha_h,0}=c_{\frac{m-h}{2},0}=\left\{\begin{array}{ll}e_h,&h\not=m,\\
e_m-{\rm dimker}\Delta,&h=m,
\end{array}\right.
\]
and we can write the expansion for the logarithmic $\Gamma$-function as 
\begin{align*}
\log\Gamma(-\lambda,Sp_+ \Delta_M)=&({\rm dimker}\Delta-e_m)\log(-\lambda)\\
&+\sum_{j=1}^{[m/2]} \frac{(-1)^{j+1}}{j!}e_{m-2j}(-\lambda)^j
\log(-\lambda)+
\sum_{h=0}^\infty a_{\scriptscriptstyle{\frac{m-h}{2}},0}
(-\lambda)^{\scriptscriptstyle{\frac{m-h}{2}}},
\end{align*}
where ($(n)$ denotes the parity of $n$)
\[
a_{\frac{m-h}{2},0}\hspace{-1pt}=\hspace{-3pt}\left\{\begin{array}{ll}\hspace{-.08in}
\Gamma\left(\frac{h-m}{2}\right)e_{h},&\hspace{-5pt}(h)\not=(m)~{\rm or}~(h)=(m)~{\rm and}~h> m,\\
\hspace{-.08in}(-1)^{\scriptscriptstyle{\frac{m-h}{2}}}
\frac{\frac{e_h}{\scriptscriptstyle{\frac{m-h}{2}}!}-\Rs{0}{\scriptscriptstyle{\frac{m-h}{2}}}\hspace{-.2in}\zeta(s,Sp_+\Delta_M)} {\scriptscriptstyle{\frac{m-h}{2}}},&
\hspace{-5pt}(h)=(m)~{\rm and}~h<m,\\
\hspace{-.08in}-\zeta'(0,Sp_+\Delta_M),&h=m,
\end{array}\right.
\]

This allows to complete the information on the zeta function as follows.

\begin{prop} The zeta function $\zeta(s,Sp_+\Delta_M)$ has an analytic continuation
to the whole complex plane up to simple poles at the values of
$s=\frac{m-h}{2}$, $h=0,1,2,\dots$, that are neither negative integers
nor zero, with residues
\[
\begin{array}{l}\Ru_{s=\frac{m-h}{2}}
\zeta(s,Sp_+ \Delta_M)=\frac{e_h}{\Gamma\left(\frac{m-h}{2}\right)}
=\left\{\begin{array}{ll}
\frac{a_{\scriptscriptstyle{\frac{h-m}{2}},0}}{\Gamma\left(\frac{h-m}{2}\right)\Gamma\left(\frac{m-h}{2}\right)},&(h)\not=(m),\\
(-1)^{\scriptscriptstyle{\frac{h-m}{2}}+1}\frac{h-m}{2}a_{\scriptscriptstyle{\frac{h-m}{2}},1}&(h)=(m),\end{array}
\right.\\
\Rz_{s=\frac{m-h}{2}}
\zeta(s,Sp_+\Delta_M)
=(-1)^{\scriptscriptstyle{\frac{h-m}{2}}+1}\frac{h-m}{2}a_{\scriptscriptstyle{\frac{h-m}{2}},1}
+\frac{e_h}{\scriptscriptstyle{\frac{m-h}{2}}!},\hspace{.2in}(h)=(m);
\end{array}
\]
the point $s=-k=0,-1,-2,\dots$ are regular points and
\begin{align*}
\zeta(0,Sp_+ \Delta_M)&=a_{0,1}=e_m-{\rm dimker}\Delta,\\
\zeta'(0,Sp_+ \Delta_M)&=-a_{0,0},\\
\zeta(-k,Sp_+ \Delta_M)&=(-1)^{k+1}ka_{k,0}=(-1)^k k! e_{m+2k}.
\end{align*}
\label{l3.1}
\end{prop}

Let $M_{(0)}=M_{(1)}\times M_{(2)}$ be the product of two compact
connected Riemannian manifolds of dimension $m_{(i)}$ without boundary. 
The metric Laplacian $\Delta_{M_{(1)}\times M_{(2)}}$ has  real spectrum with positive part
$Sp_+ \Delta_{M_{(0)}}={\{\lambda_{n_1,n_2}=\lambda_{(1),n_1}+\lambda_{(2),n_2}\}'}_{n_i\in
\N}$. Let proceed as in Section \ref{s3}. By Proposition \ref{s3.l1}, both $S_{(i)}$ are totally regular sequences of spectral type. They have exponents $s_{(i)}=\frac{m_{(i)}}{2}$, genus $p_{(i)}=\left[\frac{m_{(i)}}{2}\right]$, and infinite orders. This implies that the hypothesis of Theorem \ref{s6.t1} are satisfied,  
and consequently $S_{(0)}$ is spectrally decomposable over, say, $S_{(1)}$, with power $\kappa=1$. It also follows from the characterization  of the length given in Theorem \ref{s6.t1} and the formulas for the coefficients $\alpha_h=\frac{m-h}{2}$, just after Proposition \ref{s3.l1}, that $\ell=m_{(0)}=m_{(1)}+m_{(2)}$.
This means that we can apply Theorem \ref{tt1} and/or Corollary \ref{formulas} in order to evaluate  $\zeta(0,S_{(0)}=Sp_+ \Delta_{M_{(1)}}+Sp_+ \Delta_{M_{(2)}})$ and $\zeta'(0,S_{(0)}=Sp_+ \Delta_{M_{(1)}}+Sp_+ \Delta_{M_{(2)}})$. 
Since:
\[
\zeta(s,Sp_+ \Delta_{M_{(0)}})=\zeta(s,Sp_+ \Delta_{M_{(2)}})+\zeta(s,Sp_+ \Delta_{M_{(2)}})+\zeta(s,S_{(0)}),
\]
the following theorem follows.

\begin{teo}\label{prod} 
\[
{\rm det}_\zeta \Delta_{M_{(1)}\times M_{(2)}}={\rm det}_\zeta \Delta_{M_{(1)}}{\rm det}_\zeta \Delta_{M_{(2)}}
\e^{-\zeta'(0,Sp_+ \Delta_{M_{(1)}}+Sp_+ \Delta_{M_{(2)}})},
\]
where $\zeta'(0,Sp_+ \Delta_{M_{(1)}}+Sp_+ \Delta_{M_{(2)}})$ is given in Theorem \ref{tt1} or in Corollary \ref{formulas}.
\end{teo}

Note that a result for the zeta determinant using pure heat kernels methods  is also possible in this case. For in the case of a product manifold, one can write the regular term in the Mellin transform of the heat function by adding and subtracting the singular part of the integrand (see for example \cite {JL2} Section 3). 
This approach gives a formula for the regularized determinant involving, in the regular part, a finite integral of some complicate function, and was used in 
in \cite{FdG}. Since some derivative of the logarithmic Gamma function $\Gamma(-\lambda,S)$ is the Mellin Laplace transform of the heat function (see the proof of Proposition 2.7 of \cite{Spr2} for details), the result for $\zeta'(0,S_{(0)})$ given in Theorem \ref{tt1} is an evaluation of the integrals appearing in  Section 3 of \cite{FdG}.

\begin{exa} \label{ex3} Consider the product $S^1_{1/y}\times M$, where $S^1_r$ is the circle of radius $r$, and $M$ is a compact connected Riemannian manifold without boundary of dimension $m$, with $Sp_+ \Delta_M=\{\lambda_k\}_{k=1}^\infty$. 
Let $S_{(1)}=\{y^2n^2\}_{n=1}^\infty$, and $S_{(2)}=Sp_+ \Delta_M$. Then, $s_{(1)}=\frac{1}{2}$, $p_{(1)}=0$, $s_{(2)}=\frac{m}{2}$, and $\zeta(s,S_{(1)})=y^{-2s}\zeta_R(2s)$. Also $Sp_+\Delta_{S^1_\frac{1}{y}\times M}={\{ y^2 n^2+\lambda_k\}'}_{k\in\N,n\in \Z}$, and $S_{(0)}=S_{(1)}+S_{(2)}=\{y^2 n^2+\lambda_k\}_{n,k=1}^\infty$. It follows that  
$\zeta(s,Sp_+\Delta_{S^1_\frac{1}{y}\times M})=\zeta(s,Sp_+\Delta_M)+\zeta(s,Sp_+\Delta_{S^1_\frac{1}{y}})+2\zeta(s,S_{(0)})$.
The sequence $S_{(0)}$ is decomposable over $S_{(1)}$ with length $\ell=m+1$. We have $h=0,1,\dots,\ell=m+1$, and $\alpha_{(1),h}=\frac{1}{2},0,-\frac{1}{2},\dots,-\frac{m}{2}$. 
We apply the formula of Corollary \ref{formulas}. Since the unique pole of $\zeta(s,S_{(1)})$ is at $s=\frac{1}{2}$,  $\zeta(s,S_{(1)})$ vanishes at negative integers, and $p_{(1)}=0$, this gives
\begin{align*}
\zeta'(0,S_{(0)})=&\zeta(0,S_{(1)})\zeta'(0,S_{(2)})\\
&-2\pi \Ru_{s=\frac{1}{2}}\zeta(s,S_{(1)})\left(\Rz_{s=-\frac{1}{2}}\zeta(s,S_{(2)})+(\gamma+\psi(-\frac{1}{2})\Ru_{s=-\frac{1}{2}}\zeta(s,S_{(2)})\right)\\
&-\log\prod_{n_2=1}^\infty \e^{-2\pi \Rs{1}{\frac{1}{2}}\zeta(s,S_{(1)})\sqrt{\lambda_{(2),n_2}}}\e^{-\zeta(0,S_{(1)})\log\lambda_{(2),n_2}-\zeta'(0,S_{(1)})}\\
&\hspace{25pt}\prod_{n_1=1}^\infty \left(1+\frac{\lambda_{(2),n_2}}{\lambda_{(1),n_1}}\right).
\end{align*}
Since $\zeta(s,S_{(1)})=y^{-2s}\zeta_R(2s)$ (compare with Remark \ref{rtorus}), this gives
\begin{align*}
\zeta'(0,S_{(0)})=&-\frac{1}{2}\zeta'(0,Sp_+\Delta_M)-\log\prod_{k=1}^\infty \left(1-\e^{-\frac{2\pi}{y}\sqrt{\lambda_k}}\right)\\
&-\frac{\pi}{y}\left(\Rz_{s=-\frac{1}{2}}\zeta(s,Sp_+\Delta_M)+2(1-\log 2)\Ru_{s=-\frac{1}{2}}\zeta(s,Sp_+\Delta_M)\right),\\
&-\log\prod_{k=1}^\infty \left(1-\e^{-\frac{2\pi}{y}\sqrt{\lambda_k}}\right),
\end{align*}
and then 
\begin{align*}
\det_\zeta &\Delta_{S^1_\frac{1}{y}\times M}\\
&=\frac{4\pi^2}{y^2}
\e^{\frac{2\pi}{y}\left(\Rs{0}{-\frac{1}{2}}\hspace{-.1in}\zeta(s,Sp_+ \Delta_M)
+(2-2\log 2)\Rs{1}{-\frac{1}{2}}\hspace{-.1in}\zeta(s,Sp_+ \Delta_M)\right)}
\prod_{k=1}^\infty \left(1-\e^{-\frac{2\pi}{y}\sqrt{\lambda_k}}\right)^2.
\end{align*}

This equation is important in theoretical physics, since it gives the quantistic partition function at finite temperature $T=\frac{y}{2\pi}$ for a scalar field in the Euclidean product space time $S^1_{1/2\pi T}\times M$  (see \cite{OS} and references therein). 
\end{exa}

\section{Appendix}
\label{app}

We give in this appendix explicit formulas for some tricky countour integrals appearing in the text (see \cite{Spr3} for proofs). Here $a$ is any real number.
\begin{equation}
\frac{1}{2\pi i}\int_{\Lambda_{\theta,-c}} \e^{-\lambda }
(-\lambda)^a d \lambda=-\frac{1}{\Gamma(-a)}, \label{b1}
\end{equation}
\begin{equation}
\frac{1}{2\pi i}\int_{\Lambda_{\theta,-c}} \e^{-\lambda }
(-\lambda)^a \log (-\lambda)d
\lambda=-\frac{\psi(-a)}{\Gamma(-a)}, \label{b2}
\end{equation}
\begin{equation}\label{b4}
\int_0^\infty t^{s-1} \frac{1}{2\pi
i}\int_{\Lambda_{\theta,c}}\frac{\e^{-\lambda
t}}{-\lambda}\frac{1}{(1-\lambda)^a}d\lambda
dt=\frac{\Gamma(s+a)}{\Gamma(a)s},
\end{equation}
\begin{equation}\label{b5}
\int_0^\infty t^{s-1} \frac{1}{2\pi
i}\int_{\Lambda_{\theta,c}}\frac{\e^{-\lambda
t}}{-\lambda}\frac{\log(1-\lambda)}{(1-\lambda)^a}d\lambda
dt=\frac{\Gamma(s+a)}{\Gamma(a)s}\left(\psi(a)-\psi(s+a)\right),
\end{equation}
\beq\label{b7}
\int_0^\infty t^{s-1}\frac{1}{2\pi i}\int_{\Lambda_{\theta,c}}\frac{\e^{-\lambda t}}{-\lambda}
\log(1+\sqrt{1-\lambda})d\lambda dt=-\frac{1}{2\sqrt{\pi}}\frac{\Gamma\left(s+\frac{1}{2}\right)}{s^2}.
\eeq

\end{document}